\documentclass{article}

\usepackage{arxiv}
\usepackage[OT2, T1]{fontenc}
\usepackage[english]{babel}
\usepackage{caption}
\usepackage{lscape} 
\usepackage{amsmath}
\usepackage{geometry}
\usepackage{amsthm}
\usepackage{adjustbox}
\usepackage{amsfonts}
\usepackage{mathtools}
\usepackage{verbatim}
\usepackage{algpseudocode,algorithm,algorithmicx}
\usepackage{mathtools}
\usepackage{listings}
\usepackage{comment}
\usepackage{verbatim, times, booktabs, dcolumn, setspace, fullpage}

\newtheorem{theorem}{Theorem}

\usepackage{graphicx}
\usepackage{wrapfig}
\usepackage{caption}
\usepackage{mathtools}
\usepackage{float}

\usepackage[english,noautotitles-r]{SASnRdisplay} 
\usepackage{url}
\usepackage{subfigure}
\lstdefinestyle{r-output}{
style = r-style,
style = r-output-user,
}

\usepackage{url}            
\usepackage{booktabs}       
\usepackage{amsfonts}       
\usepackage{nicefrac}       
\usepackage{microtype}      
\usepackage{graphicx}
\usepackage[numbers]{natbib}
\usepackage{doi}

\theoremstyle{remark}
\newtheorem{remark}{Remark}

\title{A Novel Test of Missing Completely at Random: \textit{U}-statistics-based Approach}

\author{ \href{https://orcid.org/0000-0002-0460-400X}{\includegraphics[scale=0.06]{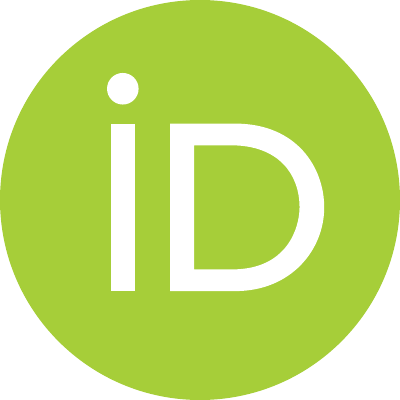}\hspace{1mm} Danijel Aleksi\' c} \\
	University of Belgrade\\
	Faculty of Organizational Sciences, Faculty of Mathematics \\
	Belgrade, 11000, Serbia \\
	\texttt{danijel.aleksic@fon.bg.ac.rs} \\
}
\date{}

\begin{document}
\maketitle

\begin{abstract}
    In this paper, a novel test for testing whether data are Missing Completely at Random is proposed. Asymptotic properties of the test are derived utilizing the theory of non-degenerate $U$-statistics. It is shown that the novel test statistic coincides with the well-known Little's $d^2$ statistic in the case of a univariate nonresponse. Then, the extensive simulation study is conducted to examine the performance of the test in terms of the preservation of type I error and in terms of power, under various underlying distributions, dimensions of the data and sample sizes. Performance of the Little's MCAR test is used as a benchmark for the comparison. The novel test shows better performance in all of the studied scenarios, better preserving the type I error and having higher empirical powers.
\end{abstract}

\section{Introduction}

Encountering missing data is a pervasive challenge in various domains, exerting an influence on the accuracy and trustworthiness of datasets. By recognizing the presence of missing data and employing effective strategies, researchers and analysts can mitigate its repercussions, ensuring more resilient and meaningful insights from the available data. The choice of appropriate techniques depends on the underlying nature of the process that resulted in the data gaps. To address different forms of missingness, it is imperative to establish a sound theoretical framework. 

First, in this section we will shortly introduce the necessary results on $U$-statistics. Then, we present the formal theoretical framework for dealing with missing data, introducing the notation that will be used throughout this paper. After that, we introduce a novel test for assessing whether the MCAR assumption holds and derive its asymptotic properties. Then, we prove that our test coincides with the well-known Little's MCAR test in the special case of univariate nonresponse. Finally, an extensive simulation study is conducted to see how our novel test performs in terms of power, using Little's MCAR test as a benchmark.

\subsection{\textit{U}-statistics}

A very important class of statistics is a class of $U$-statistics, introduced by Hoeffding in \cite{Hoeffding1948}, that turned out to be a very useful generalization of sample averages. 

Let $\mathbf{X}_1, \dots, \mathbf{X}_n$ be a sample of IID $p$-dimensional random vectors, and let  $\phi (\mathbf{x}_1, \dots, \mathbf{x}_m)$ be a measurable function  symmetric in the sense of the order of its arguments, and assume that $\mathbf{E}\phi^2(\mathbf{X}_1, \dots, \mathbf{X}_m) <\infty$. A $U$-statistic having kernel $\phi$ is defined as
\begin{align*} 
U_n = \frac{1}{\binom{n}{m}} \sum_{1 \leq i_1 < i_2 < \cdots < i_m \leq n} \phi (\mathbf{X}_{i_1}, \dots, \mathbf{X}_{i_m}). \end{align*}
If we define $\sigma_1^2$ as
\begin{align*} \sigma_1^2 = \mathbf{Cov}\left( \phi(\mathbf{X}_1, \mathbf{X}_2, \dots, \mathbf{X}_m), \phi(\mathbf{X}_1, \mathbf{X}_2', \dots, \mathbf{X}_m')  \right),\end{align*}
where $\mathbf{X}_j'$ is an independent copy of $\mathbf{X}_j$ for every $j$, we say that $U_n$ (or its kernel) is non-degenerate when $\sigma_1^2 > 0$. $U$-statistic is an unbiased and consistent estimator of $\mathbf{E}\phi(\mathbf{X}_1,\ldots,\mathbf{X}_m)$. The following theorem is a famous result proved by Hoeffding.

\begin{theorem}[Hoeffding]\label{U_statistics_joint_distribution}
Let 
\begin{align*}
    U_n^{(1)} = \frac{1}{\binom{n}{a}}\sum_{1 \leq i_1 < i_2 < \cdots < i_a \leq n} \phi^{(1)} (X_{i_1}, \dots, X_{i_a})
\end{align*}
and
\begin{align*}
    U_n = \frac{1}{\binom{n}{b}}\sum_{1 \leq i_1 < i_2 < \cdots < i_b \leq n} \phi^{(2)} (X_{i_1}, \dots, X_{i_b})
\end{align*}
be two non-degenerate $U$-statistics with kernels $\phi^{(1)}$ and $\phi^{(2)}$, respectively.
Let $\mathbf{E}\phi^{(1)} (X_1, \dots, X_a)^2 < + \infty$ and $\mathbf{E}\phi^{(2)} (X_1, \dots, X_b)^2 < + \infty$. Then 
\[ \left( \sqrt{n} \left( U_n^{(1)} - \theta_1 \right), \sqrt{n} \left( U_n - \theta_2 \right)  \right)  \overset{d}{\rightarrow} \mathcal{N}(0, \Sigma),\]
where $\Sigma$ is a limit value of the covariance matrix of $\sqrt{n} \left( U_n^{(1)} - \theta_1 \right)$ and $\sqrt{n} \left( U_n - \theta_2 \right)$, which is equal to
\begin{align*} 
\Sigma = \begin{bmatrix}
a^2 \sigma_1^{(1)} & ab \sigma_{11} \\
ab \sigma_{11} & b^2 \sigma_1^{(2)}
\end{bmatrix}, 
\end{align*}
where 
\begin{align*}
    \sigma_{11} = \mathbf{Cov}\left( \phi^{(1)} (X_1, X_2, \dots, X_a), \phi^{(2)} (X_1, X_2', \dots, X_b') \right), 
\end{align*}
 $\sigma_1^{(1)} = \mathbf{Cov}\left( \phi^{(1)} (X_1, X_2, \dots, X_{a-1}, X_a), \phi^{(1)} (X_1', X_2', \dots, X_{a-1}', X_a)  \right)$, and $\sigma_1^{(2)}$ is defined in a similar manner. Furthermore, $\Sigma$ is positive definite. 
\end{theorem}

\subsection{Theoretical framework for missing data}

Let us suppose that our data can be modeled by a $p$-dimensional random vector $\mathbf{Y} = \left( Y^{(1)}, \dots, Y^{(p)} \right)$. Now, introduce random vector $\mathbf{R} = \left( R^{(1)}, \dots, R^{(p)} \right)$, where $R^{(j)}$ is equal to $1$ if $Y^{(j)}$ is observed, and $0$ otherwise. Let us take into account an arbitrary response pattern (realisation of $\mathbf{R}$) represented by $\mathbf{r}$. We will refer to the components of $\mathbf{Y}$ that align with the positions of ones in $\mathbf{r}$ as $\mathbf{Y}_{\mathrm{obs}}$, and those aligning with zeros in $\mathbf{r}$ as $\mathbf{Y}_{\mathrm{mis}}$.

Data are said to be \emph{missing completely at random} (MCAR) if
\begin{align*} \mathbf{P} \left( \mathbf{R} = \mathbf{r} \mid  \mathbf{Y}_{\mathrm{obs}}, \mathbf{Y}_{\mathrm{mis}}\right) = \mathbf{P} \left( \mathbf{R} = \mathbf{r} \right), \end{align*}
i.e. the distribution of $\mathbf{R}$ depends neither on observed, nor on missing data, but only of its own distribution parameters. On the other hand, if the probability of observing $\mathbf{R}=\mathbf{r}$, denoted as $\mathbf{P} \left( \mathbf{R} = \mathbf{r} \mid \mathbf{Y}_{\mathrm{obs}}, \mathbf{Y}_{\mathrm{mis}}\right)$, only depends on the observed data $\mathbf{Y}_{\mathrm{obs}}$ and not on the missing data, we refer to the data as \emph{missing at random} (MAR). If either MCAR or MAR condition does not hold, we say that data are \emph{missing not at random} (MNAR). Original definitions were given by Little and Rubin in \cite{RubinLittle1987}.

One of the first missing data patterns to receive attention was a univariate missing data. That is a case where only one variable is susceptible to missingness, while the others are completely observed. This pattern is very common in many areas of research, particularly in surveys where only one question is sensitive to the participant. Well know example could be where one may refuse to disclose their income in the household survey or where patient refuses to answer to a sensitive question to a medical professional (see. e.g. \cite{RubinLittle1987}, sec. 1.6). Having been first recognized in surveys, this missing data pattern became known as univariate \emph{nonresponse}. We will need this missingness type to serve as a step to the general case of our test, which is the case where missing data can affect multiple variables.

\section{A novel test}\label{sec:ANovelTest}

The objective of this section is to present the novel test statistic and establish its asymptotic properties. We will proceed step by step, ensuring that the rationale behind the test is introduced in the most natural manner possible. First, we will cover the case of two-dimensional data with univariate nonresponse, and then gradually shift focus to the general case.

\subsection{Two-dimensional data with univaiate nonresponse}

Let us suppose that we have data from two dimensional distribution that can be represented by the random vector $(X, Y)$, where $\mathbf{E}(X^2) < \infty$, and let us model our sample as its $n$ IID copies, which we expand to obtain the sample
\begin{align}\label{expanded_sample}  \begin{bmatrix}
X_1 & Y_1 & R_1 \\
X_2 & Y_2 & R_2 \\
\vdots & \vdots & \vdots \\
X_n & Y_n & R_n
\end{bmatrix},
\end{align}
where every $X_j$ is observed and $R_j$ denotes the response indicator for $Y_j$, i.e.
\begin{align*} 
R_j = \begin{cases}
    1, & \text{ if } Y_j \text{ is observed,} \\
    0, & \text{otherwise.}
\end{cases} 
\end{align*}

If data are missing completely at random (MCAR), then, by definition, it holds that $\mathbf{P} (R \mid X, Y_{obs}, Y_{miss}) = \mathbf{P} (R)$. But, in such case, it trivially holds that $\mathbf{P} (X \mid R) = \mathbf{P} (X)$, i.e. $\mathbf{P} (XR) = \mathbf{P} (X) \mathbf{P} (R)$, which means that $X$ and $R$ are independent. As a consequence, under MCAR assumption it holds that $\mathbf{E}(XR) = \mathbf{E}(X) \mathbf{E}(R)$, or, equivalently,
\begin{align*}  
\mathbf{E}(X) \mathbf{E}(R) - \mathbf{E}(XR) = 0.
\end{align*}
Naturally, for the distribution of $(X, Y, R)$ we can define parameter $\theta = \mathbf{E}(X) \mathbf{E}(R) - \mathbf{E}(XR)$; if $\theta \neq 0$, then $X$ and $R$ are dependent, so the data are not MCAR. If $\theta = 0$, than we can not say that the independence holds, but only that $X$ and $R$ are not correlated. However, there are many well-known tests that reject the null hypothesis in such manner (e.g. Kendall's test of independence, see e.g. \cite{Kendall1970}). This being said, it is natural to construct test for MCAR based on some estimator of a parameter $\theta$. One such estimator is 
\begin{align}\label{Tn_tilde}
    \Tilde{T}_n = \left(\frac{1}{n} \sum_i X_i \right)  \left(  \frac{1}{n} \sum_j R_j \right) - \frac{1}{n} \sum_i X_iR_i = \frac{1}{n^2} {\sum_{i=1}^n \sum_{\substack{j = 1 \\ j \neq i}}^n} X_iR_j - \frac{n-1}{n^2}\sum_{i=1}^n  X_i R_i.
\end{align}
$\Tilde{T}_n$ is a biased estimator. Indeed, one can easily verify that $\mathbf{E}(\Tilde{T}_n) = \frac{n-1}{n} \theta$. After appropriate rescaling, we obtain an unbiased estimator of $\theta$ given by:
\begin{align}\label{Tn}
    T_n = \frac{1}{n(n-1)} {\sum_{i=1}^n \sum_{\substack{j = 1 \\ j \neq i}}^n} X_iR_j - \frac{1}{n} \sum_{i=1}^n X_iR_i.
\end{align}
After going through some transformations, we have:
\begin{align*}
    T_n &= \frac{1}{\binom{n}{2}} \underset{1 \leq i < j \leq n}{\sum \sum} \frac12 (X_iR_j + X_j R_i) - \frac{1}{\binom{n}{1}} \sum_{i=1}^n X_iR_i.
\end{align*}
If we denote
\begin{align*}
    U_n^{(1)} = \frac{1}{\binom{n}{2}} \underset{1 \leq i < j \leq n}{\sum \sum} \phi \left( (X_i, Y_i, R_i), (X_j, Y_j, R_j)  \right)
\end{align*}
and
\begin{align*}
    U_n^{(2)} = \frac{1}{\binom{n}{1}} \sum_{i=1}^n \psi((X_i, Y_i, R_i)),
\end{align*}
where $\phi \left( (X_i, Y_i, R_i), (X_j, Y_j, R_j)  \right) = \frac12 (X_iR_j + X_jR_i)$ and $\psi ((X_i, Y_i, R_i)) = X_i R_i$, we have that
\begin{align*}
    T_n = U_n^{(1)} - U_n^{(2)}.
\end{align*}
Note that $U_n^{(1)}$ and $U_n^{(2)}$ are $U$-statistics that have kernels $\phi$ and $\psi$, respectively. The next step will be to show that they are non-degenerate, and to obtain asymptotic distribution of $T_n$ under the null hypothesis.

We can see that
\begin{align*}
    \sigma_1^{(2)} &= \mathbf{Cov}(\psi (X_1, Y_1, R_1), \psi(X_1, Y_1, R_1)) = \mathbf{Var}\left(\psi (X_1, Y_1, R_1)\right) =  \mathbf{E}(X^2 R^2) - \mathbf{E}(X)^2 \mathbf{E}(R)^2.
\end{align*}
Under the null hypothesis, we have that $X$ and $R$ are independent variables, and after we note that $R^2 = R$, we obtain 
\begin{align*}
    \sigma_1^{(2)} = \mathbf{E}(X^2) \mathbf{E}(R) - \mathbf{E}(X)^2 \mathbf{E}(R)^2.
\end{align*}
Since we are interested in the non-trivial case (not all data observed/missing), we can safely assume that $\mathbf{E}(R) \in (0,1)$, so $\mathbf{E}(R) > \left( \mathbf{E}(R) \right)^2$. Furthermore, $\mathbf{E}(X^2)  > \left(\mathbf{E}(X) \right)^2$. These two conditions allow us to conclude that $\sigma_1^{(2)} > 0$, i.e. $U_n^{(2)}$ is a non-degenerate $U$-statistic.

Now we proceed in a similar manner with statistic $U_n^{(1)}$ and, after similar calculations as before, obtain:
\begin{align*}
    \sigma_1^{(1)} &= \mathbf{Cov}\left( \phi ((X_1, Y_1, R_1), (X_2, Y_2, R_2)),  \phi ((X_1, Y_1, R_1), (X_3, Y_3, R_3))   \right) \\
    &= \frac14 \left(\mathbf{E}(R)^2 \mathbf{Var}(X) + \mathbf{E}(X)^2 \mathbf{Var}(R) \right).
\end{align*}
Again, being interested only in non-trivial cases, we have that $\sigma_1^{(1)} > 0$, so $U_n^{(1)}$ is a non-degenerate $U$-statistic.

It is not difficult for one to see that $U_n^{(1)}$ is a $U$-statistic that estimates the parameter $\mathbf{E}\left( \phi \left( (X_1, Y_1, R_1), (X_2, Y_2, R_2)  \right) \right) = \mathbf{E}(X) \mathbf{E}(R)$, and that $U_n^{(2)}$ is a $U$-statistic that estimates the parameter $\mathbf{E}\left( \psi(X_1, Y_1, R_1) \right) = \mathbf{E}(XR)$. By the known formula for covariance of the two $U$-statistics (see, e.g. \cite{Hoeffding1948}), we have that
\begin{align*}
    \mathbf{Cov}\left(U_n^{(1)}, U_n^{(2)}\right) = \frac{1}{\binom{n}{1}} \sum_{i=1}^1 \binom{2}{i} \binom{n-2}{2-i} \sigma_{ii} = \frac{2(n-2)}{n} \sigma_{11},
\end{align*}
where $\sigma_{11} = \mathbf{Cov}(\phi((X_1, Y_1, R_1), (X_2, Y_2, R_2)), \psi (X_1, Y_1, R_1))$. Under the null hypothesis we have that:
\begin{align*}
    \sigma_{11} &= \mathbf{Cov}\left( \frac12 ( X_1R_2 + X_2R_1), X_1R_1 \right) \\
    &= \frac12 \left(\mathbf{E}\left( (X_1R_2 + X_2R_1)X_1R_1  \right) - \mathbf{E}(X_1R_2 + X_2R_1) \mathbf{E}(X_1R_1) \right) \\
    &= \frac12 \left( \mathbf{E}(R)^2 \mathbf{Var}(X) + \mathbf{E}(X)^2 \mathbf{Var}(R) \right).
\end{align*}

Having Theorem \ref{U_statistics_joint_distribution}, we have that as $n \to \infty$, under the null hypothesis,
\begin{align}\label{asimptotika_scriptsize}
      \bigg(  \sqrt{n} \big(& U_n^{(1)} - \mathbf{E}(X)\mathbf{E}(R)\big),    \sqrt{n} \big(U_n^{(2)} - \mathbf{E}(XR)\big)  \bigg) \\ \nonumber
    & \overset{D}{\to}  \mathcal{N}_2 \left( \mathbf{0}, \begin{bmatrix}
     \mathbf{E}(R)^2 \mathbf{Var}(X) + \mathbf{E}(X)^2 \mathbf{Var}(R)  &  \mathbf{E}(R)^2 \mathbf{Var}(X) + \mathbf{E}(X)^2 \mathbf{Var}(R) \\
    \mathbf{E}(R)^2 \mathbf{Var}(X) + \mathbf{E}(X)^2 \mathbf{Var}(R) & \mathbf{E}(X^2) \mathbf{E}(R) - \mathbf{E}(X)^2 \mathbf{E}(R)^2
\end{bmatrix}     \right). 
\end{align}
Applying Continuous mapping theorem, we have that difference $\sqrt{n} (U_n^{(1)} - U_n^{(2)})$ (expectations cancel out under null hypothesis) tends to difference of the components of the two-dimensional normal distribution given in \eqref{asimptotika_scriptsize}, i.e. 
\begin{align*}
    \sqrt{n} \big(U_n^{(1)} - U_n^{(2)}  \big) \\ 
    \overset{D}{\to} \mathcal{N} \bigg(  0,\, &   \mathbf{E}(R)^2 \mathbf{Var}(X) + \mathbf{E}(X)^2 \mathbf{Var}(R)  \\
    & \quad\quad - 2\big(\mathbf{E}(R)^2 \mathbf{Var}(X) + \mathbf{E}(X)^2 \mathbf{Var}(R)\big) + \mathbf{E}(X^2) \mathbf{E}(R) - \mathbf{E}(X)^2 \mathbf{E}(R)^2  \bigg), 
\end{align*}
which, after some canceling, becomes
\begin{align*}
    \sqrt{n} \left(U_n^{(1)} - U_n^{(2)}  \right)  \overset{D}{\to} \mathcal{N} \left(  0,   \mathbf{Var}(X) \mathbf{Var}(R) \right). 
\end{align*}
In other words,
\begin{align*}
    \frac{\sqrt{n} \left(U_n^{(1)} - U_n^{(2)}  \right)}{\sqrt{\mathbf{Var}(X) \mathbf{Var}(R)}} \overset{D}{\to} \mathcal{N} (0,1),
\end{align*}
as $n \to \infty$. Knowing that (adjusted for bias) sample standard deviation is a consistent estimator, then applying Slutsky's theorem, we have that
\begin{align*}
    \frac{\sqrt{n} \left(U_n^{(1)} - U_n^{(2)}  \right)}{S_n^{X} S_n^R} \overset{D}{\to} \mathcal{N} (0,1).
\end{align*}
Having this, we suggest constructing a MCAR test using test statistic
\begin{align*}
    D_n = \frac{\sqrt{n} \left(U_n^{(1)} - U_n^{(2)}  \right)}{S_n^{X} S_n^R},
\end{align*}
and having rejection region
\begin{align}\label{Dn}
     \bigg\{ |D_n| \geq \Phi^{-1} \left( 1 - \frac{\alpha}{2}  \right) \bigg\}
\end{align}
for a significance level $\alpha$, where $\Phi$ is a cumulative distribution function of the standard normal distribution.

\subsection{Multivariate data with univariate nonresponse}

Now we will go a step further from the two-dimensional data, and  consider $p$-dimensional data with univariate nonresponse. To keep the expressions simple, we will consider $p = 3$, but it will be obvious that the generalization to the arbitrary $p$ is straightforward. 

Consider the data that can be modeled by random vector $(X^{(1)}, X^{(2)}, Y)$ and, as before, consider the expanded sample
\begin{align}\label{expanded_sample_2}  
\begin{bmatrix}
X^{(1)}_1 & X^{(2)}_1 & Y_1 & R_1 \\
X^{(1)}_2 & X^{(2)}_2 & Y_2 & R_2 \\
\vdots & \vdots & \vdots & \vdots \\
X^{(1)}_n & X^{(2)}_n & Y_n & R_n
\end{bmatrix}.
\end{align}
At this point, we assume that only $Y$ is susceptible to missingness. It is natural to utilize statistic $T_n$ for the variable pairs $(X^{(1)}, Y)$ and $(X^{(2)}, Y)$. Remembering \eqref{Tn}, this consideration leads us to the statistics
\begin{align}\label{Tn1}
    T_n^{(1)} = \frac{1}{n(n-1)} {\sum_{i=1}^n \sum_{\substack{j = 1 \\ j \neq i}}^n} X_i^{(1)}R_j - \frac{1}{n} \sum_{i=1}^n X_i^{(1)}R_i 
\end{align}
and
\begin{align}\label{Tn2}
    T_n^{(2)} = \frac{1}{n(n-1)} {\sum_{i=1}^n \sum_{\substack{j = 1 \\ j \neq i}}^n} X_i^{(2)}R_j - \frac{1}{n} \sum_{i=1}^n X_i^{(2)}R_i. 
\end{align}
In the previous section, we have proven that both statistics, scaled by $\sqrt{n}$, have asymptotic normal distribution, being differences of the non-degenerate $U$-statistics. Small values of either $T_n^{(1)}$ or $T_n^{(2)}$ indicate that the null hypothesis of MCAR is true, so, combining the two, it is natural to associate the null hypothesis being whenever both $T_n^{(1)}$ and $T_n^{(2)}$ are small.

Wanting to construct the test based on this fact, we  first obtain the asymptotic distribution of the random vector $\left( \sqrt{n} T_n^{(1)}, \sqrt{n} T_n^{(2)} \right)$. Certainly, it would be wrong, in the general case, to assume joint normality from the normality of the components. However, since the statistics $\sqrt{n} T_n^{(1)}$ and $ \sqrt{n}T_n^{(2)}$ are (scaled) differences of $U$-statistics, their joint distribution can be written as a linear combination of joint distributions of $U$-statistics, and hence is asymptotically normal. Since we now have that the asymptotic distribution of $\left( \sqrt{n} T_n^{(1)}, \sqrt{n} T_n^{(2)} \right)$ is normal, it will suffice to calculate the limit value of the covariance
\begin{align*}
    \mathbf{Cov}\left( \sqrt{n}T_n^{(1)}, \sqrt{n} T_n^{(2)} \right) = n \mathbf{Cov}\left( T_n^{(1)}, T_n^{(2)} \right) = n \mathbf{E}\left( T_n^{(1)} T_n^{(2)} \right),
\end{align*}
as $n \to \infty$.

For the calculations, we utilize similar technique found in \cite{aleksic2023etAl}. Multiplying the expressions \eqref{Tn1} and \eqref{Tn2}, we obtain
\begin{align*}
    T_n^{(1)} T_n^{(2)} &= \frac{1}{n^2 (n-1)^2} \sum_{i=1}^n \sum_{\substack{j = 1 \\ j \neq i}}^n \sum_{k = 1}^n \sum_{\substack{l = 1 \\ l \neq k}}^n X_i^{(1)} R_j X_k^{(2)}R_l + \frac{1}{n^2} \sum_{i=1}^n \sum_{j=1}^n X_i^{(1)} R_i X_j^{(2)} R_j \\
    & \quad \quad- \frac{1}{n^2(n-1)} \sum_{i=1}^n \sum_{\substack{j=1\\ j \neq i}}^n \sum_{k=1}^n X_i^{(1)} R_j X_k^{(2)} R_k - \frac{1}{n^2(n-1)} \sum_{k=1}^n \sum_{\substack{l = 1 \\ l \neq k}}^n \sum_{i=1}^n X_k^{(2)}R_l X_i^{(1)} R_i \\
    &=: M + N - Q_1 - Q_2.
\end{align*}
Noting that
\begin{align*}
    n^2(n-1)^2 M &= \sum_i \sum_{j \neq i} \sum_{\substack{l \neq i \\ l \neq j}} X_i^{(1)} R_j X_i^{(2)} R_l + \sum_{i} \sum_{j \neq i} X_i^{(1)} R_j X_i^{(2)} \\
    &+ \sum_i \sum_{j \neq i} \sum_{\substack{k \neq i \\ k \neq j}} \sum_{\substack{l \neq i \\ l \neq j \\ l \neq k}} X_i^{(1)} R_j X_k^{(2)} R_l + \sum_{i} \sum_{j \neq i} \sum_{\substack{k \neq i \\ k \neq j}} X_i^{(1)} R_j X_k^{(2)} R_i + \sum_i \sum_{j \neq i} \sum_{\substack{k \neq i \\ k \neq j}} X_i^{(1)} R_j X_k^{(2)} \\
    &+ \sum_i \sum_{j \neq i} \sum_{\substack{l \neq i \\ l \neq j}} X_i^{(1)} R_j X_j^{(2)} R_l + \sum_i \sum_{j \neq i} X_i^{(1)} R_j X_j^{(2)} R_i,
\end{align*}
we can obtain
\begin{align*}
    \mathbf{E}\left(  n^2(n-1)^2 M \right) &= n(n-1)(n-2) \mathbf{E}\left( X^{(1)} X^{(2)} \right) \left( \mathbf{E}(R) \right)^2 \\
    &+ n(n-1) \mathbf{E}\left(  X^{(1)} X^{(2)} \right) \mathbf{E}(R) \\
    & + n(n-1)(n-2)(n-3) \mathbf{E}\left( X^{(1)} \right) \mathbf{E}\left( X^{(2)} \right) \left( \mathbf{E}(R) \right)^2 \\
    &+ n(n-1)(n-2) \mathbf{E}\left( X^{(1)} \right) \mathbf{E}\left( X^{(2)} \right) \left( \mathbf{E}(R) \right)^2 \\
    &+ n(n-1)(n-2) \mathbf{E}\left( X^{(1)} \right) \mathbf{E}\left( X^{(2)} \right) \mathbf{E}(R) \\
    &+ n(n-1)(n-2) \mathbf{E}\left( X^{(1)} \right) \mathbf{E}\left( X^{(2)} \right) \left( \mathbf{E}(R) \right)^2 \\
    &+ n(n-1) \mathbf{E}\left( X^{(1)} \right) \mathbf{E}\left( X^{(2)} \right) \left( \mathbf{E}(R) \right)^2.
\end{align*}
Having this, we can see that
\begin{align*}
    \mathbf{E}\left( nM \right) &= \frac{n-2}{n-1} \mathbf{E}\left( X^{(1)} X^{(2)} \right) \left( \mathbf{E}(R) \right)^2 + \frac{(n-2)(n-3)}{n-1}\mathbf{E}\left( X^{(1)} \right) \mathbf{E}\left( X^{(2)} \right) \left( \mathbf{E}(R) \right)^2 \\
    &+ 2 \frac{n-2}{n-1} \mathbf{E}\left( X^{(1)} \right) \mathbf{E}\left( X^{(2)} \right) \left( \mathbf{E}(R) \right)^2 + \frac{n-2}{n-1} \mathbf{E}\left( X^{(1)} \right) \mathbf{E}\left( X^{(2)} \right) \left( \mathbf{E}(R) \right) + o(1), 
\end{align*}
as $n \to \infty$. In a similar manner, it follows that
\begin{align*}
    \mathbf{E}\left( nN \right) &= (n-1) \mathbf{E}\left( X^{(1)} \right) \mathbf{E}\left( X^{(2)} \right) \left( \mathbf{E}(R) \right)^2 + \mathbf{E}\left( X^{(1)} X^{(2)} \right) \mathbf{E}(R)
\end{align*}
and 
\begin{align*}
    \mathbf{E}(nQ_1) = \mathbf{E}(nQ_2) &= 2(n-2)\mathbf{E}\left( X^{(1)} \right) \mathbf{E}\left( X^{(2)} \right) \left( \mathbf{E}(R) \right)^2 + 2 \mathbf{E}\left( X^{(1)} X^{(2)} \right) \left(\mathbf{E}(R) \right)^2 \\
    & \quad \quad + 2\mathbf{E}\left( X^{(1)} \right) \mathbf{E}\left( X^{(2)} \right) \left( \mathbf{E}(R) \right).
\end{align*}
Combining previous results, and since $\mathbf{E}\left( T_n^{(1)} \right) = \mathbf{E}\left( T_n^{(2)} \right) = 0$ and $R^2 = R$, we have that
\begin{align*}
    \lim_{n \to \infty} \mathbf{Cov}\left( \sqrt{n} T_n^{(1)}, \sqrt{n} T_n^{(2)} \right) &= \lim_{ n \to \infty} n \mathbf{E}\left( T_n^{(1)} T_n^{(2)} \right) \\
    &= \lim_{n \to \infty} \left( nM + nN - nQ1 - nQ2 \right) \\
    &= \mathbf{E}\left( X^{(1)}  X^{(2)} \right) \left( \mathbf{E}(R) \right) - \mathbf{E}\left( X^{(1)}  X^{(2)} \right) \mathbf{E}(R)  \\
    &+ \mathbf{E}\left( X^{(1)} \right) \mathbf{E}\left( X^{(2)} \right) \left( \mathbf{E}(R) \right)^2 - \mathbf{E}\left( X^{(1)}  X^{(2)} \right) \left( \mathbf{E}(R) \right)^2 \\
    &= \mathbf{Cov}\left( X^{(1)}, X^{(2)} \right) \mathbf{E}\left( R^2 \right) - \mathbf{Cov}\left( X^{(1)}, X^{(2)} \right) \left( \mathbf{E}\left( R \right) \right)^2 \\
    &= \mathbf{Cov}\left( X^{(1)}, X^{(2)} \right)\mathbf{Var}(R).
\end{align*}
At this point, we have successfully proven the following theorem.
\begin{theorem}\label{Teorema3d}
    Let us have data represented by the expanded sample \eqref{expanded_sample_2}, and let $T_n^{(1)}$ and $T_n^{(2)}$ be given in \eqref{Tn1} and \eqref{Tn2}, respectively. Then, it holds that 
    \begin{align*}
        \left( \sqrt{n} T_n^{(1)}, \sqrt{n} T_n^{(2)} \right) \overset{D}{\to} \mathcal{N} \left( \boldsymbol{0}, \begin{bmatrix}
            \mathbf{Var}\left( X^{(1)} \right) \mathbf{Var}(R) & \mathbf{Cov}\left( X^{(1)}, X^{(2)} \right)\mathbf{Var}(R) \\
            \mathbf{Cov}\left( X^{(1)}, X^{(2)} \right)\mathbf{Var}(R) & \mathbf{Var}\left( X^{(2)} \right) \mathbf{Var}(R)
        \end{bmatrix} \right),
    \end{align*}
    as $n \to \infty$.
\end{theorem}

Now, suppose that we have the $(p+1)$-variate data, that we observe as the expanded sample
\begin{align}\label{expanded_sample_p}
    \begin{bmatrix}
X^{(1)}_1 & X^{(2)}_1 & \cdots & X^{(p)}_1 & Y_1 & R_1 \\
X^{(1)}_2 & X^{(2)}_2 & \cdots & X^{(p)}_2 & Y_2 & R_2 \\
\vdots & \vdots & \ddots & \vdots & \vdots & \vdots \\
X^{(1)}_n & X^{(2)}_n & \cdots & X^{(p)}_n & Y_n & R_n
\end{bmatrix}.
\end{align}
The next theorem is the straightforward generalization of the Theorem \ref{Teorema3d}. 
\begin{theorem}\label{TeoremaPpd}
     Let us have data represented by the expanded sample \eqref{expanded_sample_p}, and let 
     \begin{align*}
         T_n^{(u)} = \frac{1}{n(n-1)} \sum_{i=1}^n \sum_{\substack{j=1 \\ j \neq i}}^n X_i^{(u)}R_j - \frac{1}{n} \sum_{i=1}^n X_i^{(u)}R_i,
     \end{align*}
     $u = 1, 2, \dots, p$. Then, it holds that
     \begin{align*}
         \left( \sqrt{n} T_n^{(1)}, \sqrt{n} T_n^{(2)}, \dots, \sqrt{n} T_n^{(p)} \right) \overset{D}{\to} \mathcal{N} \left( \boldsymbol{0}, \boldsymbol{\Sigma}  \right),
     \end{align*}
     as $n \to \infty$, where
     \begin{align*}
         \boldsymbol{\Sigma} = \begin{bmatrix}
            \mathbf{Var}\left( X^{(1)} \right) \mathbf{Var}(R) & \mathbf{Cov}\left( X^{(1)}, X^{(2)} \right)\mathbf{Var}(R) & \cdots &  \mathbf{Cov}\left( X^{(1)}, X^{(p)} \right)\mathbf{Var}(R)\\
            \mathbf{Cov}\left( X^{(1)}, X^{(2)} \right)\mathbf{Var}(R) & \mathbf{Var}\left( X^{(2)} \right) \mathbf{Var}(R) & \cdots & \mathbf{Cov}\left( X^{(p)}, X^{(2)} \right)\mathbf{Var}(R) \\
            \vdots & \vdots & \ddots & \vdots \\
            \mathbf{Cov}\left( X^{(1)}, X^{(p)} \right)\mathbf{Var}(R) & \mathbf{Cov}\left( X^{(2)}, X^{(p)} \right)\mathbf{Var}(R) \cdots &\cdots & \mathbf{Var}\left( X^{(p)} \right) \mathbf{Var}(R)
        \end{bmatrix},
     \end{align*}
     i.e.
     \begin{align*}
         \boldsymbol{\Sigma} = \mathbf{Cov}\left( X^{(1)}, X^{(2)}, \dots, X^{(p)} \right) \mathbf{Var}(R),
     \end{align*}
     where first term denotes the covariance matrix of a random vector.
\end{theorem}

Now, let us introduce another assumption: there is no such pair $\left( X^{(u)}, X^{(v)} \right)$ in the vector $\left(X^{(1)}, X^{(2)}, \dots, X^{(p)} \right)$ such that $\mathbf{Cov}\left( X^{(u)}, X^{(v)} \right) /\sqrt{\mathbf{Var}(X^{(u)}) \mathbf{Var}(X^{(v)})} = \pm 1$, for $u, v = 1, 2, \dots, p$, $u \neq v$. It is a well-known result (see e.g. \cite{strang2016introduction}, p. 549) that in that case $\boldsymbol{\Sigma}$ is positive definite matrix, and that there exists the matrix $\boldsymbol{\Sigma}^{-\frac{1}{2}}$ such that $\boldsymbol{\Sigma}\cdot \left( \boldsymbol{\Sigma}^{-\frac{1}{2}}\right)^2 = \left( \boldsymbol{\Sigma}^{-\frac{1}{2}}\right)^2 \cdot  \boldsymbol{\Sigma} = \mathbf{E}$, where $\mathbf{E}$ is the identity matrix. Since matrix multiplication is linear, and hence continuous transformation, we can use Continuous mapping theorem to obtain 
\begin{align*}
    \left(\boldsymbol{\Sigma}^{-\frac{1}{2}} \cdot \left( \sqrt{n} T_n^{(1)}, \sqrt{n} T_n^{(2)}, \dots, \sqrt{n} T_n^{(p)} \right)^T \right)^T \overset{D}{\to} \mathcal{N} (\boldsymbol{0}, \mathbf{E}),
\end{align*}
as $n \to \infty$.

Assume additionally that every variable $X^{(u)}$, $u = 1,2,\dots, p$ has finite fourth moment. Components of the matrix $\boldsymbol{\Sigma}$ are generally unknown and need to be estimated. That can be done using standard bias-adjusted sample covariance matrix, multiplied by bias-adjusted sample variance of $R$, that are known to be consistent estimators whenever fourth moments of the variables are finite. That leads us to estimated matrix $\boldsymbol{\hat{\Sigma}}$. Since inverse and square root of a matrix are obtained by linear, and hence continuous transformations, by the Continuous mapping theorem it holds that $\boldsymbol{\hat{\Sigma}}^{-\frac12}$ is a consistent estimator of $\boldsymbol{\Sigma}^{-\frac12}$. Using Slutsky's theorem, applied to every component, we obtain that
\begin{align}\label{raspodela_Ap}
    \left( A_n^{(1)}, A_n^{(2)}, \dots, A_n^{(p)} \right) := \left(\boldsymbol{\hat{\Sigma}}^{-\frac{1}{2}} \cdot \left( \sqrt{n} T_n^{(1)}, \sqrt{n} T_n^{(2)}, \dots, \sqrt{n} T_n^{(p)} \right)^T \right)^T\overset{D}{\to} \mathcal{N} (\boldsymbol{0}, \mathbf{E}),
\end{align}
as $n \to \infty$.

Values of any $T_n^{(u)}$, $u = 1, 2, \dots, p$, close to zero are a strong indicator of a MCAR assumption being true. The same holds for the components of a vector $\left( A_n^{(1)}, A_n^{(2)}, \dots, A_n^{(p)} \right)$, it being just a linear transformation of the previous. Having this conclusion, we construct a test based on a statistic
\begin{align*}
    A_n = \left(A_n^{(1)}\right)^2 + \left(A_n^{(2)}\right)^2 + \cdots + \left(A_n^{(p)}\right)^2,
\end{align*}
whose small values indicate that MCAR holds. Having \eqref{raspodela_Ap}, we see that $A_n$ is asymptotically distributed as a sum of squares of $p$ IID standard normal variables, and hence $A_n$ is asymptotically $\chi_p^2$-distributed. We construct a critical region of a test by rejecting the null hypothesis if $A_n$ is larger than the adequate quantile of a $\chi_p^2$ distribution.

\begin{remark}\label{remark_sigma}
    Note that the statistic $A_n$, being equal to sum of the squares of the vector 
    \begin{align*}
    \boldsymbol{\hat{\Sigma}}^{-\frac12} \cdot \left( \sqrt{n} T_n^{(1)}, \dots, \sqrt{n}T_n^{(p)} \right)^T,
    \end{align*}
    can be written as
    \begin{align*}
        A_n = n \left( T_n^{(1)}, T_n^{(2)}, \dots, T_n^{(p)} \right) \boldsymbol{\hat{\Sigma}}^{-1} \left( T_n^{(1)}, T_n^{(2)}, \dots, T_n^{(p)} \right)^T.
    \end{align*}
    If we, instead, used maximum likelihood estimate $\boldsymbol{\Tilde{\Sigma}}$ of a matrix $\boldsymbol{\Sigma}$, that is not unbiased, we would have relation $\boldsymbol{\hat{\Sigma}} = \left(\frac{n}{n-1}\right)^2 \boldsymbol{\Tilde{\Sigma}}$, and we could write
    \begin{align}\label{An_mat}
        A_n = n \left( \tilde{T}_n^{(1)}, \tilde{T}_n^{(2)}, \dots, \tilde{T}_n^{(p)} \right) \boldsymbol{\Tilde{\Sigma}}^{-1} \left( \tilde{T}_n^{(1)}, \tilde{T}_n^{(2)}, \dots, \tilde{T}_n^{(p)} \right)^T,
    \end{align}
    where $\Tilde{T}_n^{(u)}$, $u = 1, 2, \dots, p$, are defined analogously as in \eqref{Tn_tilde}.
\end{remark}

\subsection{The general case}

As a final step, we will soften the assumption of a univariate nonresponse, and allow multiple variables to be susceptible to missing data. The only restriction we impose is that there is at least one completely observed variable. Consider the data that can be modeled by random vector $(X^{(1)}, \dots, X^{(p)}, Y^{(1)}, \dots, Y^{(q)})$ and consider the expanded sample
\begin{align}\label{expanded_sample_last}  
\begin{bmatrix}
X^{(1)}_1 & \cdots & X^{(p)}_1 & Y_1^{(1)} & \cdots & Y_1^{(q)} & R_1^{(1)} & \cdots & R_1^{(q)} \\
X^{(1)}_2 & \cdots & X^{(p)}_2 & Y_2^{(1)} & \cdots & Y_2^{(q)} & R_2^{(1)} & \cdots & R_2^{(q)} \\
\vdots & \ddots &\vdots & \vdots & \ddots & \vdots & \vdots & \vdots \\
X^{(1)}_n & \cdots & X^{(p)}_n & Y_n^{(1)} & \cdots  & Y_n^{(q)}  & R_n^{(1)} & \cdots & R_n^{(q)}
\end{bmatrix}.
\end{align}
Suppose that variables $X^{(1)}, \dots, X^{(p)}$ are completely observed and that $Y^{(1)}, \dots, Y^{(q)}$ are susceptible to missingness.  Similarly to previous subsections, we introduce the statistics 
\begin{align}\label{Tn_last}
    T_n^{(u, v)} =  \frac{1}{n(n-1)} {\sum_{i=1}^n \sum_{\substack{j = 1 \\ j \neq i}}^n} X_i^{(u)}R_j^{(v)} - \frac{1}{n} \sum_{i=1}^n X_i^{(u)}R_i^{(v)},
\end{align}
$u = 1, 2, \dots, p$, $v = 1, 2, \dots, q$. In the exact same manner as the Theorem \ref{TeoremaPpd}, we obtain the following theorem:
\begin{theorem}\label{theorem_final}
Let us have data represented by the expanded sample \eqref{expanded_sample_last} and let $T_n^{(u,v)}$ be defined as in \eqref{Tn_last}, for $u = 1, 2, \dots, p$ and $v = 1,2, \dots, q$. Then, it holds that
\begin{align*}
\left(  \sqrt{n} T_n^{(1,1)}, \dots, \sqrt{n} T_n^{(1, q)}, \sqrt{n} T_n^{(2,1)}, \dots, \sqrt{n} T_n^{(2,q)}, \dots, \sqrt{n}T_n^{(p,1)}, \dots, \sqrt{n}T_n^{(p,q)}  \right) \overset{D}{\to} \mathcal{N} (\boldsymbol{0}, \boldsymbol{\Sigma}),
\end{align*}
as $n \to \infty$, where
\begin{align*}
    \boldsymbol{\Sigma} = \left[ \mathbf{Cov}\left( X^{(\lceil i / p \rceil)} X^{(\lceil j / p \rceil)} \right), \mathbf{Cov}\left( R^{(i \, (\mathrm{mod}\, p))}, R^{(j \, (\mathrm{mod}\, p))} \right) \right]_{\substack{i,j \in \{1, \dots, pq\}}}
\end{align*}
where $a\, \mathrm{mod}\, b$ is remainder of the division of $a$ by $b$, and $\lceil \cdot \rceil$ is the ceiling function. 
\end{theorem}
As before, if we assume that complete variables have finite fourth moments, $\boldsymbol{\Sigma}^{-\frac12}$ exists and can be consistently estimated in a similar manner to obtain $\boldsymbol{\hat{\Sigma}}^{-\frac{1}{2}}$ or $\boldsymbol{\Tilde{\Sigma}}^{-\frac{1}{2}}$ as a consistent estimator. Finally, it holds that
\begin{align*}
    \big( A_n^{(1,1)}, \dots, A_n^{(1, q)}, \dots, A_n^{(2,1)}, \dots,& A_n^{(2, q)}, \dots, A_n^{(p,1)}, \dots, A_n^{(p, q)}  \big) \\
    &= \left(\boldsymbol{\hat{\Sigma}}^{-\frac12} \cdot \left( \sqrt{n} T_n^{(1,1)}, \dots, \sqrt{n}T_n^{(p,q)} \right)^T\right)^T \overset{D}{\to} \mathcal{N} (\boldsymbol{0}, \mathbf{E}),
\end{align*}
as $n \to \infty$, so we have that statistic
\begin{align*}
    A_n = \sum_{u = 1}^p \sum_{v=1}^q \left( A_n^{(u ,v)} \right)^2 
\end{align*}
asymptotically has $\chi^2_{pq}$ distribution, that can be used to calculate the critical values for the test.

\begin{remark}
    Note that as long we have $p$ complete and $q$ incomplete variables, the data can be rearranged to take the form \eqref{expanded_sample_last}.
\end{remark}

\section{A note on the special case of univariate nonresponse}

One of the most well-known tests for testing the MCAR assumption is Little's MCAR test, constructed by Little in \cite{Little1988}. For a sample of independent and identically distributed random vectors, with univariate nonresponse, as in \eqref{expanded_sample_p}, we shall prove that ours and Little's statistics coincide. For that purpose, we first need to adapt the expression for Little's statistic to our notation. For brevity, we discuss the three-dimensional case \eqref{expanded_sample_2}, but it will be obvious that the generalization is straightforward.

Denote the vectors
\begin{align*}
    L = \left( \frac{1}{\sum_{i=1}^n R_i}\sum_{i=1}^n X_i^{(1)} R_i - \frac{1}{n} \sum_{i=1}^n X_i^{(1)},  
    \frac{1}{\sum_{i=1}^n R_i} \sum_{i=1}^n X_i^{(2)} R_i - \frac{1}{n} \sum_{i=1}^n X_i^{(2)}\right),
\end{align*}
\begin{align*}
    L_1 = \left( \frac{1}{\sum_{i=1}^n (1 - R_i)}\sum_{i=1}^n X_i^{(1)} (1 - R_i) - \frac{1}{n} \sum_{i=1}^n X_i^{(1)},  
     \frac{1}{\sum_{i=1}^n (1 - R_i)}\sum_{i=1}^n X_i^{(2)} (1 - R_i) - \frac{1}{n} \sum_{i=1}^n X_i^{(2)}, 0\right)
\end{align*}
and
\begin{align*}
    L_1' = \left( \frac{1}{\sum_{i=1}^n (1 - R_i)}\sum_{i=1}^n X_i^{(1)} (1 - R_i) - \frac{1}{n} \sum_{i=1}^n X_i^{(1)},  
     \frac{1}{\sum_{i=1}^n (1 - R_i)}\sum_{i=1}^n X_i^{(2)} (1 - R_i) - \frac{1}{n} \sum_{i=1}^n X_i^{(2)}\right),
\end{align*}
and let $\boldsymbol{\Tilde{\Sigma}}_1$ be a matrix obtained by expanding $\boldsymbol{\Tilde{\Sigma}}$ to include estimates of covariance of $X^{(i)}$, $i = 1, 2$, and $Y$, calculated on those rows $i$ that have $Y_i$ observed. So, the matrix $\boldsymbol{\Tilde{\Sigma}}_1$ is of a form 
\begin{align*}
    \boldsymbol{\Tilde{\Sigma}}_1 = \begin{bmatrix}
        \boldsymbol{\Tilde{\Sigma}} & \boldsymbol{\Tilde{\Lambda}} \\
        \boldsymbol{\Tilde{\Lambda}}^T & \boldsymbol{\Tilde{\Delta}}
    \end{bmatrix},
\end{align*}
where $\boldsymbol{\Tilde{\Lambda}}$ and $\boldsymbol{\Tilde{\Delta}}$ are estimators of adequate matrices, and $\boldsymbol{\Tilde{\Sigma}}$ is defined in Remark \ref{remark_sigma}, here just for $p = 2$. Denote $\boldsymbol{\Tilde{\Sigma}}' = \frac{1}{\Bar{R}_n (1 - \Bar{R}_n)} \boldsymbol{\Tilde{\Sigma}}$ and $\boldsymbol{\Tilde{\Sigma}}_1' = \frac{1}{\Bar{R}_n (1 - \Bar{R}_n)} \boldsymbol{\Tilde{\Sigma}}_1$, where $\Bar{R}_n = \frac{1}{n} \sum_{i=1}^n R_i$.

In his paper \cite{Little1988}, Little introduced the test statistic that, in this special case of univariate nonresponse, has a form
\begin{align*}
    d^2 = n \Bar{R}_n L (\boldsymbol{\Tilde{\Sigma}}')^{-1} L^T + n (1 - \Bar{R}_n) L_1 (\boldsymbol{\Tilde{\Sigma}}_1')^{-1} L_1^T.
\end{align*}
But, since in this case vector $L_1$ has zero as the last component, matrices $\boldsymbol{\Tilde{\Lambda}}$ and $\boldsymbol{\Tilde{\Delta}}$ have no effect on the $d^2$, so we can make reduction and obtain that
\begin{align*}
    d^2 = n \Bar{R}_n L (\boldsymbol{\Tilde{\Sigma}}')^{-1} L^T + n (1 - \Bar{R}_n) L_1' (\boldsymbol{\Tilde{\Sigma}}')^{-1} L_1'^T,
\end{align*}
which, after substituting $\boldsymbol{\Tilde{\Sigma}}'$ and $\boldsymbol{\Tilde{\Sigma}}_1'$, become
\begin{align}\label{d2}
    d^2 = n \Bar{R}_n^2 (1 - \Bar{R}_n) L \boldsymbol{\Tilde{\Sigma}}^{-1} L^T + n \Bar{R}_n (1 - \Bar{R}_n)^2 L_1' \boldsymbol{\Tilde{\Sigma}}^{-1} L_1'^T.
\end{align}

Denote $\Bar{X}_n^{(1)} = \frac{1}{n} \sum_{i=1}^n X_i^{(1)}$ and $\Bar{X}_n^{(2)}$ similarly. Also, let $\overline{X_n^{(1)}R_n} = \frac{1}{n} \sum_{i=1}^n X_i^{(1)}R_i$, and let $\overline{X_n^{(2)}R_n}$ be defined in an analogous manner. Now we can see that 
\begin{align*}
    L =  \frac{1}{\Bar{R}_n} \left( \overline{X_n^{(1)} R_n}, \overline{X_n^{(2)} R_n}  \right) - \left( \Bar{X}_n^{(1)}, \Bar{X}_n^{(2)} \right)
\end{align*}
and 
\begin{align*}
    L_1' &= \frac{1}{(1 - \Bar{R}_n)} \left(  \Bar{X}_n^{(1)} - \overline{X_n^{(1)} R_n}, \Bar{X}_n^{(2)} - \overline{X_n^{(2)} R_n} \right) - \left( \Bar{X}_n^{(1)}, \Bar{X}_n^{(2)} \right) \\
    &= \frac{\Bar{R}_n}{1 - \Bar{R}_n} \left( \Bar{X}_n^{(1)}, \Bar{X}_n^{(2)} \right) - \frac{1}{(1 - \Bar{R}_n)} \left( \overline{X_n^{(1)} R_n}, \overline{X_n^{(2)} R_n}  \right).
\end{align*}
Next, we calculate both of the terms in \eqref{d2}. First, we have that
\begin{align*}
   \Bar{R}_n^2 (1 - \Bar{R}_n) L \boldsymbol{\Tilde{\Sigma}}^{-1} L^T &= n (1 - \Bar{R}_n) \left( \overline{X_n^{(1)} R_n}, \overline{X_n^{(2)} R_n}  \right) \boldsymbol{\Tilde{\Sigma}}^{-1} \left( \overline{X_n^{(1)} R_n}, \overline{X_n^{(2)} R_n}  \right)^T \\
    & \quad \quad \quad -2 n \Bar{R}_n (1 - \Bar{R}_n) \left( \overline{X_n^{(1)} R_n}, \overline{X_n^{(2)} R_n}  \right) \boldsymbol{\Tilde{\Sigma}}^{-1}\left( \Bar{X}_n^{(1)}, \Bar{X}_n^{(2)} \right)^T \\
    & \quad \quad \quad + n \Bar{R}_n^2 (1 - \Bar{R}_n) \left( \Bar{X}_n^{(1)}, \Bar{X}_n^{(2)} \right) \boldsymbol{\Tilde{\Sigma}}^{-1}\left( \Bar{X}_n^{(1)}, \Bar{X}_n^{(2)} \right)^T.
\end{align*}
Then, it holds that
\begin{align*}
    n \Bar{R}_n (1 - \Bar{R}_n)^2 L_1' \boldsymbol{\Tilde{\Sigma}}^{-1} L_1'^T &= n \Bar{R}_n \left( \overline{X_n^{(1)} R_n}, \overline{X_n^{(2)} R_n}  \right) \boldsymbol{\Tilde{\Sigma}}^{-1} \left( \overline{X_n^{(1)} R_n}, \overline{X_n^{(2)} R_n}  \right)^T \\
    & \quad \quad \quad - 2 n \Bar{R}_n^2 \left( \overline{X_n^{(1)} R_n}, \overline{X_n^{(2)} R_n}  \right) \boldsymbol{\Tilde{\Sigma}}^{-1}\left( \Bar{X}_n^{(1)}, \Bar{X}_n^{(2)} \right)^T \\
    & \quad \quad \quad + n \Bar{R}_n^3 \left( \Bar{X}_n^{(1)}, \Bar{X}_n^{(2)} \right)\boldsymbol{\Tilde{\Sigma}}^{-1} \left( \Bar{X}_n^{(1)}, \Bar{X}_n^{(2)} \right)^T.
\end{align*}
Combining, we have
\begin{align*}
    d^2 &= n \left( \overline{X_n^{(1)} R_n}, \overline{X_n^{(2)} R_n}  \right) \boldsymbol{\Tilde{\Sigma}}^{-1} \left( \overline{X_n^{(1)} R_n}, \overline{X_n^{(2)} R_n}  \right)^T \\
    & \quad \quad \quad -2 n \Bar{R}_n \left( \overline{X_n^{(1)} R_n}, \overline{X_n^{(2)} R_n}  \right) \boldsymbol{\Tilde{\Sigma}}^{-1}\left( \Bar{X}_n^{(1)}, \Bar{X}_n^{(2)} \right)^T \\
    & \quad \quad \quad + n \Bar{R}_n^2 \left( \Bar{X}_n^{(1)}, \Bar{X}_n^{(2)} \right)\boldsymbol{\Tilde{\Sigma}}^{-1} \left( \Bar{X}_n^{(1)}, \Bar{X}_n^{(2)} \right)^T.
\end{align*}

On the other hand, from \eqref{An_mat} (for $p = 2$) and from the fact that $\Tilde{T}_n^{(u)} = \Bar{X}_n^{(u)} \Bar{R}_n - \overline{X_n^{(u)} R_n}$, $u = 1, 2$, it follows that 
\begin{align*}
    A_n &= n \left( \tilde{T}_n^{(1)}, \tilde{T}_n^{(2)} \right) \boldsymbol{\Tilde{\Sigma}}^{-1} \left( \tilde{T}_n^{(1)}, \tilde{T}_n^{(2)} \right)^T\\
    &= n \left( \overline{X_n^{(1)} R_n}, \overline{X_n^{(2)} R_n}  \right) \boldsymbol{\Tilde{\Sigma}}^{-1} \left( \overline{X_n^{(1)} R_n}, \overline{X_n^{(2)} R_n}  \right)^T \\
    & \quad \quad \quad -2 n \Bar{R}_n \left( \overline{X_n^{(1)} R_n}, \overline{X_n^{(2)} R_n}  \right) \boldsymbol{\Tilde{\Sigma}}^{-1}\left( \Bar{X}_n^{(1)}, \Bar{X}_n^{(2)} \right)^T \\
    & \quad \quad \quad + n \Bar{R}_n^2 \left( \Bar{X}_n^{(1)}, \Bar{X}_n^{(2)} \right)\boldsymbol{\Tilde{\Sigma}}^{-1} \left( \Bar{X}_n^{(1)}, \Bar{X}_n^{(2)} \right)^T,
\end{align*}
which is exactly $d^2$. Thus, the following theorem holds.

\begin{theorem}\label{theorem_equal}
    For the $(p+1)$-variate data with univariate nonresponse, that can be represented as expanded sample \eqref{expanded_sample_p}, it holds that
    \begin{align*}
        A_n = d^2,
    \end{align*}
    where $A_n$ is as in \eqref{An_mat}, and $d^2$ is defined in \eqref{d2}.
\end{theorem}

\section{Empirical study}

After we have derived the needed statistical properties of our novel test, we proceed to verify that test behaves well in terms of size preservation and power. For that purpose, a simulation study is conducted to examine the behavior. Then our test is compared to the Little's MCAR test, using it as a benchmark.

We generate \textbf{datasets of several dimensions} and missingness structure. Namely, the three different ones are used. First, we generate three-dimensional data with one complete and two incomplete variables, which we will refer to as $1X2Y$ \emph{case}. Then, five-dimensional data is generated, with three complete and two incomplete variables, which we refer to as $3X2Y$ \emph{case}. Last, the five-dimensional data with two complete and three incomplete columns is studied - the $2X3Y$ \emph{case}. 

Besides the dimensions and patterns, we go across the set of \textbf{different underlying distributions}. First, the multivariate standard normal distribution as a point of reference, knowing that Little's test has some of it's properties derived for the normal case. Then, following reasoning from \cite{CuparicMilosevic2023}, to explore different inner data structures, Clayton copula with parameter $1$ (which gives the Kendall's correlation between variables of $1/3$) is used, both with standard exponential margins and with chi-squared margins with $4$ degrees of freedom. 

\textbf{MCAR data} are generated such that every value has the same probability to be missing, independent of others. We use implementation from the \texttt{R} package \texttt{missMethods}. For the size comparison, we first compare empirical sizes as functions of missingness probability, for fixed dimension, distribution, and sample size. The probability that a value is missing is ranged from $0.03$ to $0.24$, with the step of $0.03$. Then, we compare empirical sizes as functions of sample size, for the fixed missingness probability, distribution and dimension, to study the speed of convergence to the asymptotic distribution.

For the alternatives to MCAR, 
the implementations of \textbf{two MAR mechanisms} are used, that were originally and formally introduced in \cite{santos2019generating}. We refer to the same \texttt{missMethods} package for the implementations. The idea behind them is natural: every incomplete column has its complete pair, which dictates the missingness. For the first mechanism, which is called \emph{MAR 1 to x}, if value it the complete column is higher than its median, than the pair value in incomplete column is $x$ times more likely to be missing. We use the suggested choice $x=9$, that is default in \texttt{R}. The second mechanism is so-called \emph{MAR rank mechanism}, where probability that a value is missing is directly dependent on the rank of its observed pair. For detailed explanations of mechanisms, one should consult \cite{santos2019generating}. Besides these two, we use \textbf{another MAR mechanism} that does not need missingness probability as a parameter. This mechanism was used by the authors of \cite{mozharovskyi2020nonparametric}, for studying their technique of imputation by data depth. The mechanism was introduced for the data of a form 1X2Y, and the missingness probability was set different for the units that have $X$ larger than its mean value and for those that have it smaller. We will refer to this mechanism as the \emph{MAR mean} mechanism.

\subsection{Size behavior}

If one wants the inference about power to be valid, the test needs to preserve its type I error. We use the standard choice of $0.05$. 

\begin{figure}[ht]
    \centering
    \includegraphics[width=\textwidth]{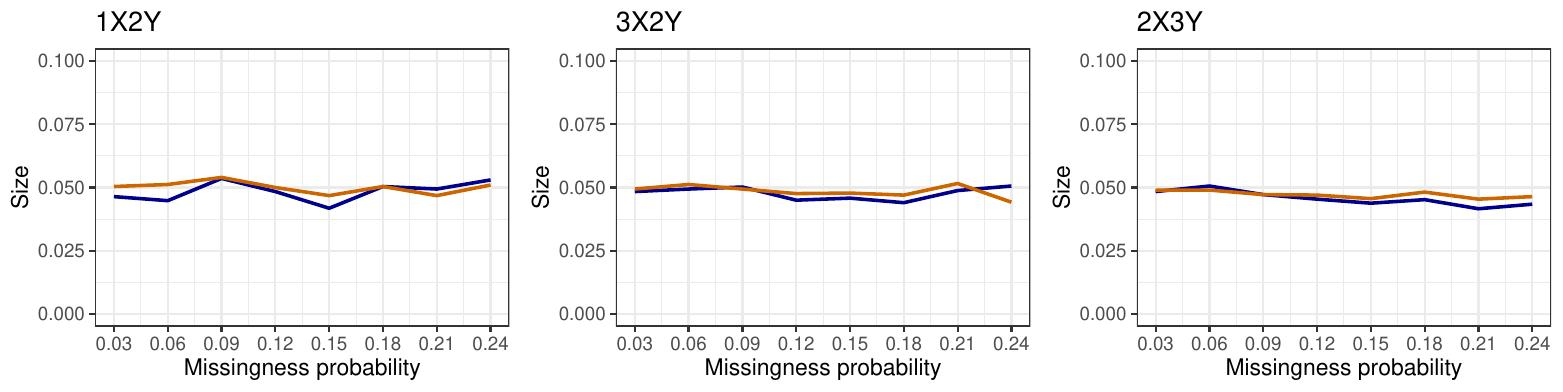}
    \caption{Empirical test \textbf{sizes} ($A_n$ in orange, $d^2$ in blue) for \textbf{standard normal} distribution, sample size $\mathbf{n = 100}$, and $N = 5000$ simulations }
    \label{fig:sizes_stdnorm_n100}
\end{figure}

\begin{figure}[ht]
    \centering
    \includegraphics[width=\textwidth]{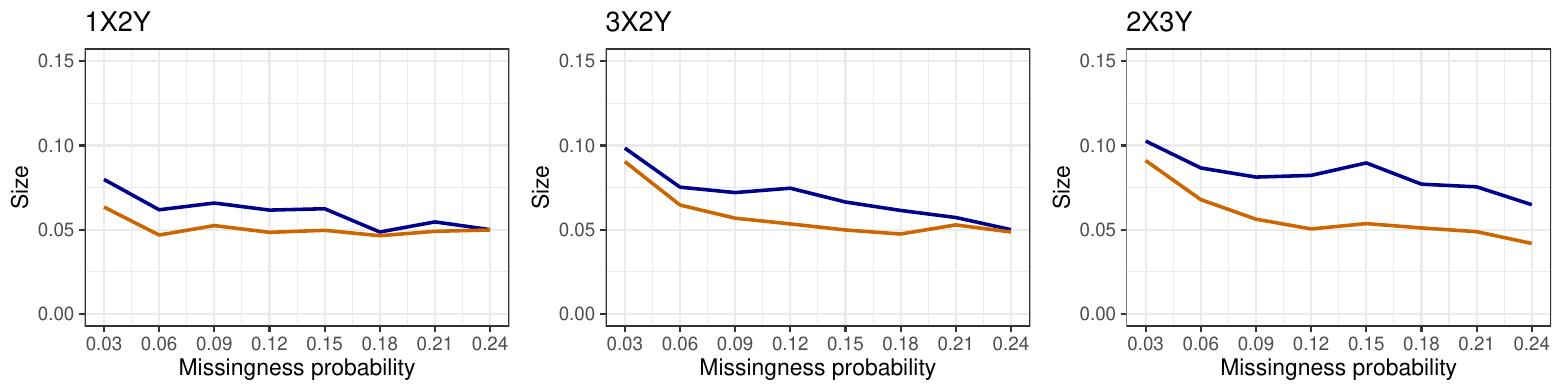}
    \caption{Empirical test \textbf{sizes} ($A_n$ in orange, $d^2$ in blue) for \textbf{Clayton copula with parameter 1} and $\mathbf{\mathcal{E} (1)}$ \textbf{margins}, sample size $\mathbf{n = 100}$, and $N = 5000$ simulations }
    \label{fig:sizes_clayton1_exp1_n100}
\end{figure}

\begin{figure}[ht]
    \centering
    \includegraphics[width=\textwidth]{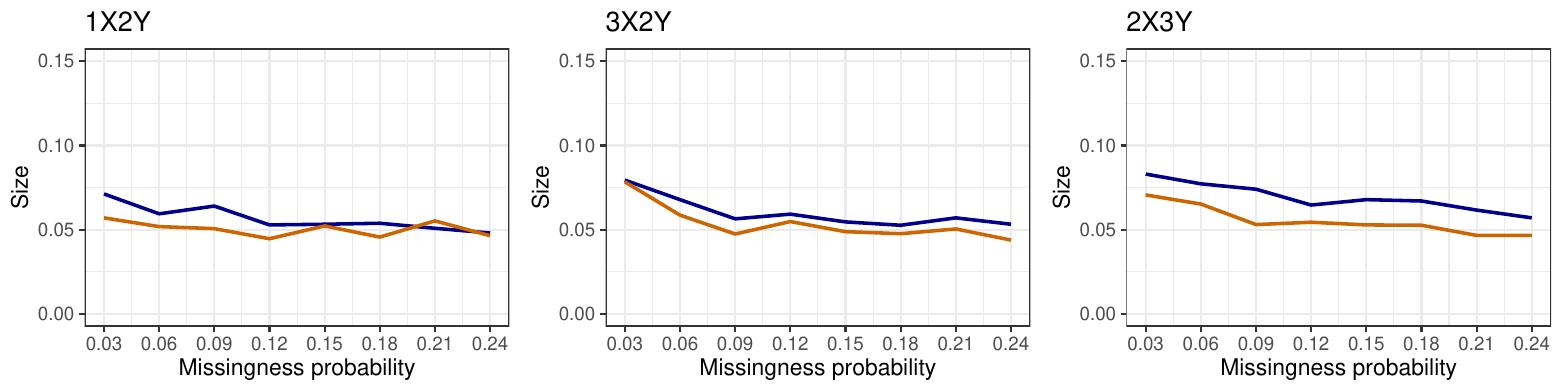}
    \caption{Empirical test \textbf{sizes} ($A_n$ in orange, $d^2$ in blue) for \textbf{Clayton copula with parameter 1} and $\mathbf{\chi^2_4}$ \textbf{margins}, sample size $\mathbf{n = 100}$, and $N = 5000$ simulations }
    \label{fig:sizes_clayton1_chisq4_n100}
\end{figure}

Figure \ref{fig:sizes_stdnorm_n100} shows us that both $A_n$ and $d^2$ preserve the test size for the sample size of 100 when underlying distribution is standard normal. However, as we can see from the Figure \ref{fig:sizes_clayton1_exp1_n100}, when the underlying distribution is not normal, sizes start to differ. Test constructed around statistic $A_n$ experiences slight empirical size increase for the small missingness probabilities, but stabilizes around the probability of $0.1$ for the worst observed case. However, Little's MCAR test experiences significant increase in size, especially for five-dimensional data with trivariate nonresponse. This behavior persists if we change marginal distributions from standard exponential to chi-squared, as can be seen in the figure \ref{fig:sizes_clayton1_chisq4_n100}. We see that $d^2$ performs a little bit better in this case, since $\chi^2$ distributions have higher level of similarity to the normal distribution than the exponential ones. 

\begin{figure}[ht]
    \centering
    \includegraphics[width=\textwidth]{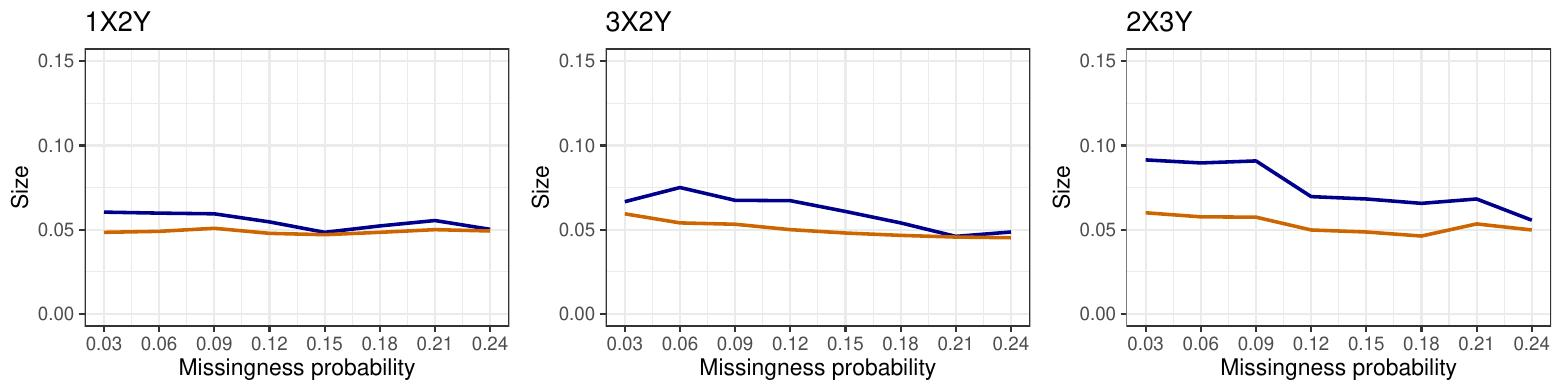}
    \caption{Empirical test \textbf{sizes} ($A_n$ in orange, $d^2$ in blue) for \textbf{Clayton copula with parameter 1} and $\mathbf{\mathcal{E} (1)}$ \textbf{margins}, sample size $\mathbf{n = 300}$, and $N = 5000$ simulations }
    \label{fig:sizes_clayton1_exp1_n300}
\end{figure}

These results tell us that asymptotic convergence of Little's $d^2$ statistic is significantly slowed down if the underlying distribution is not standard normal. The inference for the test power against alternatives would not be valid in this case, so it is necessary for a sample size to be increased. As can be seen from the Figure \ref{fig:sizes_clayton1_exp1_n300}, the problem persists, although a little bit less notable, for the sample size $ n = 300$ and Clayton copula with exponential margins. 

\begin{figure}[ht]
    \centering
    \includegraphics[width=\textwidth]{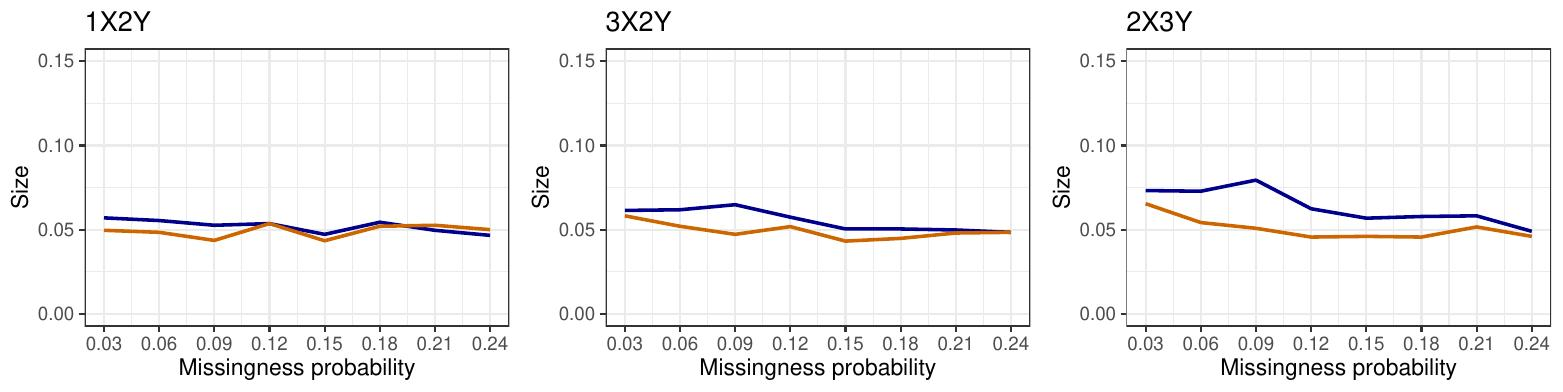}
    \caption{Empirical test \textbf{sizes} ($A_n$ in orange, $d^2$ in blue) for \textbf{Clayton copula with parameter 1} and $\mathbf{\chi^2_4}$ \textbf{margins}, sample size $\mathbf{n = 300}$, and $N = 5000$ simulations }
    \label{fig:sizes_clayton1_chisq4_n300}
\end{figure}

Figure \ref{fig:sizes_clayton1_chisq4_n300} shows us that the problem with statistic $d^2$ is tolerable for chi-squared margins and sample size of $n = 300$. We note that several other cases of normal distributions and Clayton copula were tried, as well as different data dimensions. However, since they behaved in the same manner as those presented here, their inclusion in the paper would be redundant.

\begin{figure}[ht]
    \centering
    \includegraphics[width=\textwidth]{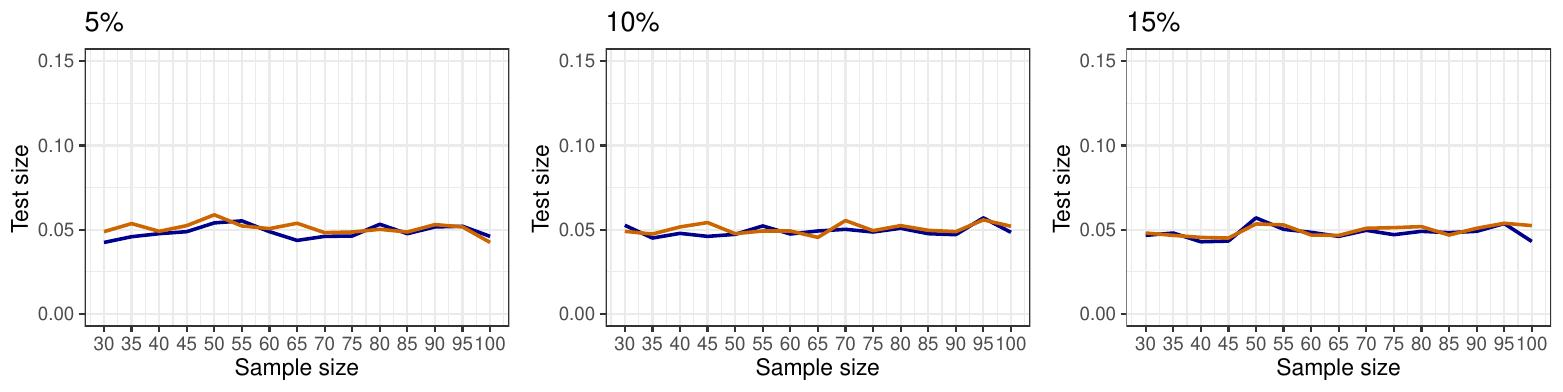}
    \caption{Empirical test \textbf{sizes for 1X2Y case} ($A_n$ in orange, $d^2$ in blue) for \textbf{standard normal} distribution, \textbf{sample size ranging from 30 to 100}, $N = 5000$ simulations }
    \label{fig:sizes_1X2Y_stdnorm_n30to100}
\end{figure}

\begin{figure}[ht]
    \centering
    \includegraphics[width=\textwidth]{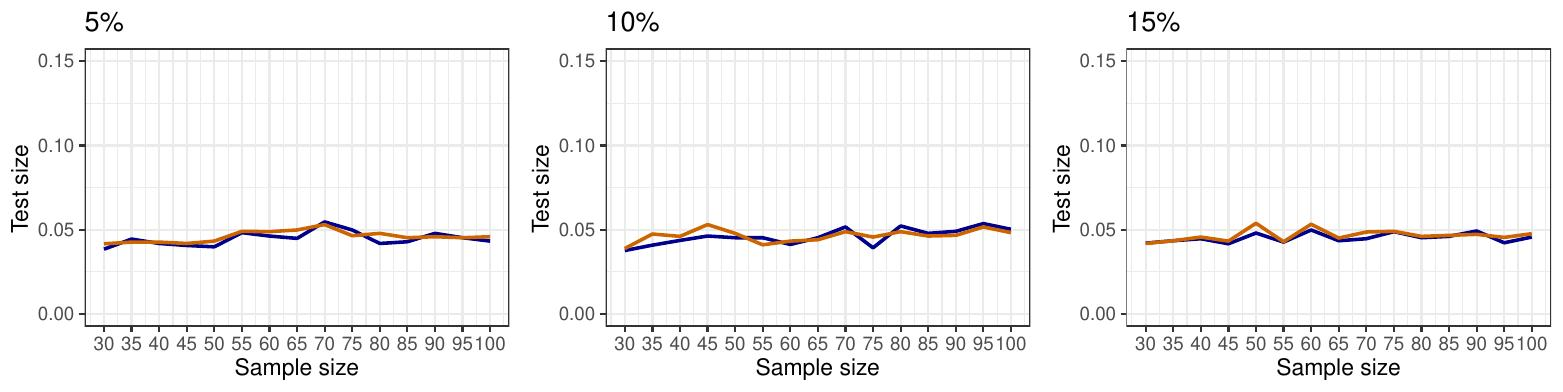}
    \caption{Empirical test \textbf{sizes for 3X2Y case} ($A_n$ in orange, $d^2$ in blue) for \textbf{standard normal} distribution, \textbf{sample size ranging from 30 to 100}, $N = 5000$ simulations }
    \label{fig:sizes_3X2Y_stdnorm_n30to100}
\end{figure}

\begin{figure}[ht]
    \centering
    \includegraphics[width=\textwidth]{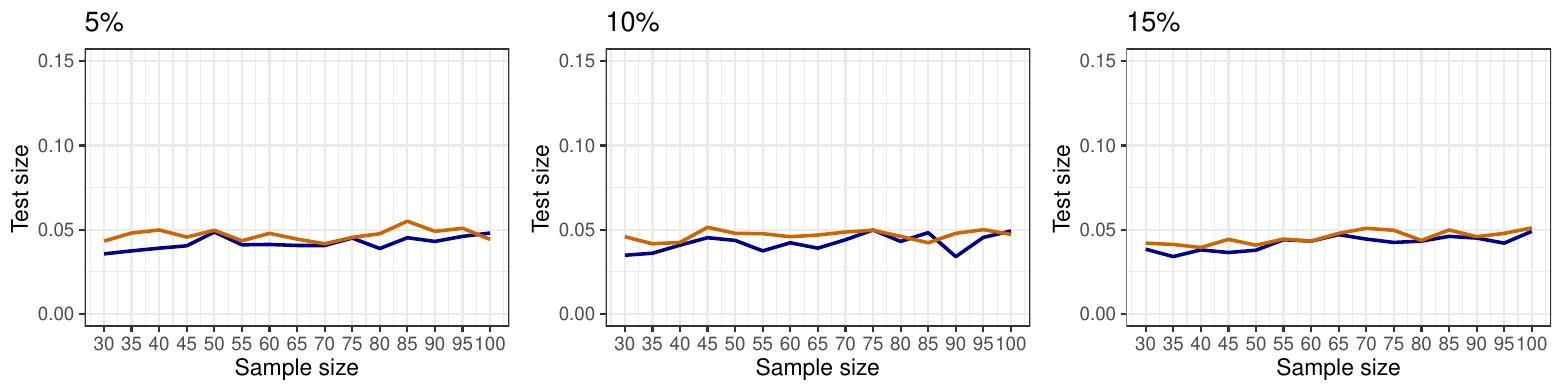}
    \caption{Empirical test \textbf{sizes for 2X3Y case} ($A_n$ in orange, $d^2$ in blue) for \textbf{standard normal} distribution, \textbf{sample size ranging from 30 to 100}, $N = 5000$ simulations }
    \label{fig:sizes_2X3Y_stdnorm_n30to100}
\end{figure}

Next, we compare the two tests in terms of the speed of convergence to the asymptotic distribution. The data distribution, dimension, and probability that a value is missing are fixed, and test size for various sample sizes, ranging from 30 to 100 is compared. As we can see from the Figures \ref{fig:sizes_1X2Y_stdnorm_n30to100}, \ref{fig:sizes_3X2Y_stdnorm_n30to100} and \ref{fig:sizes_2X3Y_stdnorm_n30to100}, both tests converge very quickly for the underlying standard normal distribution - size is preserved even for the sample as small as 30 observations. Figure \ref{fig:sizes_2X3Y_stdnorm_n30to100} shows us that $A_n$ in the 2X3Y case has slightly faster convergence than $d^2$, which in this case underestimates the size, but not significantly.

\begin{figure}[ht]
    \centering
    \includegraphics[width=\textwidth]{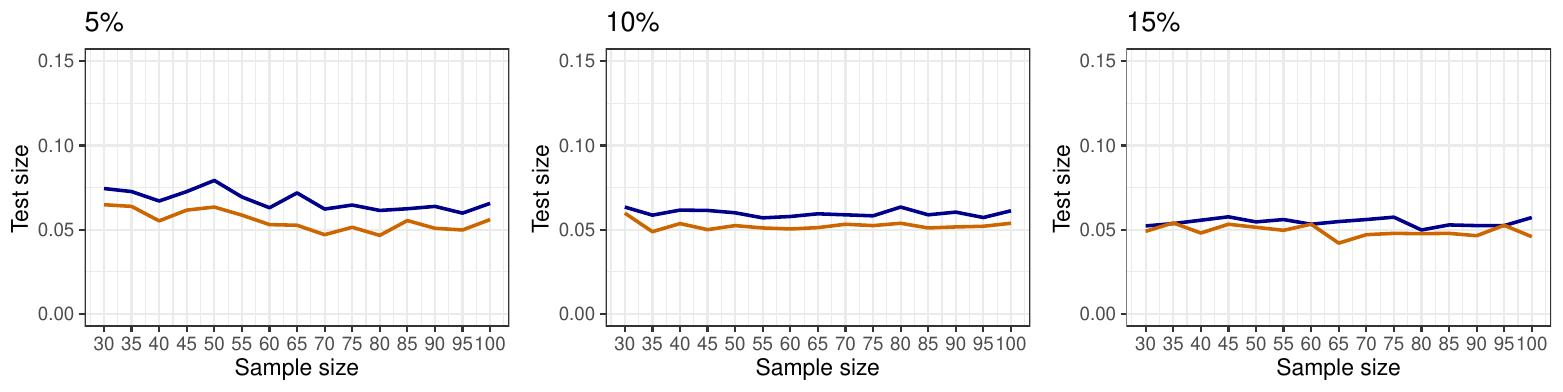}
    \caption{Empirical test \textbf{sizes for 1X2Y case} ($A_n$ in orange, $d^2$ in blue) for \textbf{Clayton copula with parameter 1} and $\mathbf{\chi^2_4}$ \textbf{margins}, \textbf{sample size ranging from 30 to 100}, $N = 5000$ simulations }
    \label{fig:sizes_1X2Y_clayton1_chisq4_n30to100}
\end{figure}

\begin{figure}[ht]
    \centering
    \includegraphics[width=\textwidth]{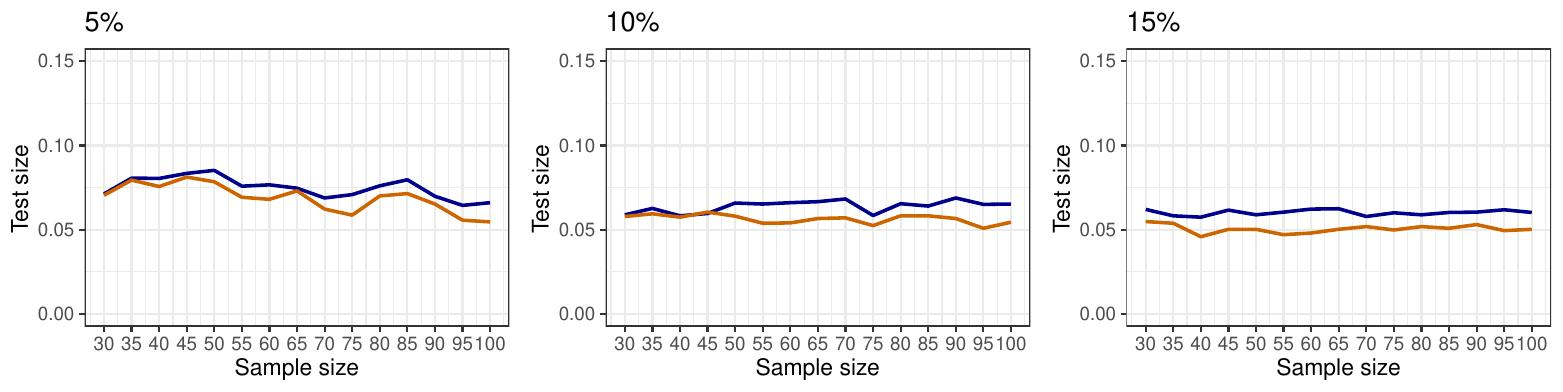}
    \caption{Empirical test \textbf{sizes for 3X2Y case} ($A_n$ in orange, $d^2$ in blue) for \textbf{Clayton copula with parameter 1} and $\mathbf{\chi^2_4}$ \textbf{margins}, \textbf{sample size ranging from 30 to 100}, $N = 5000$ simulations }
    \label{fig:sizes_3X2Y_clayton1_chisq4_n30to100}
\end{figure}

\begin{figure}[ht]
    \centering
    \includegraphics[width=\textwidth]{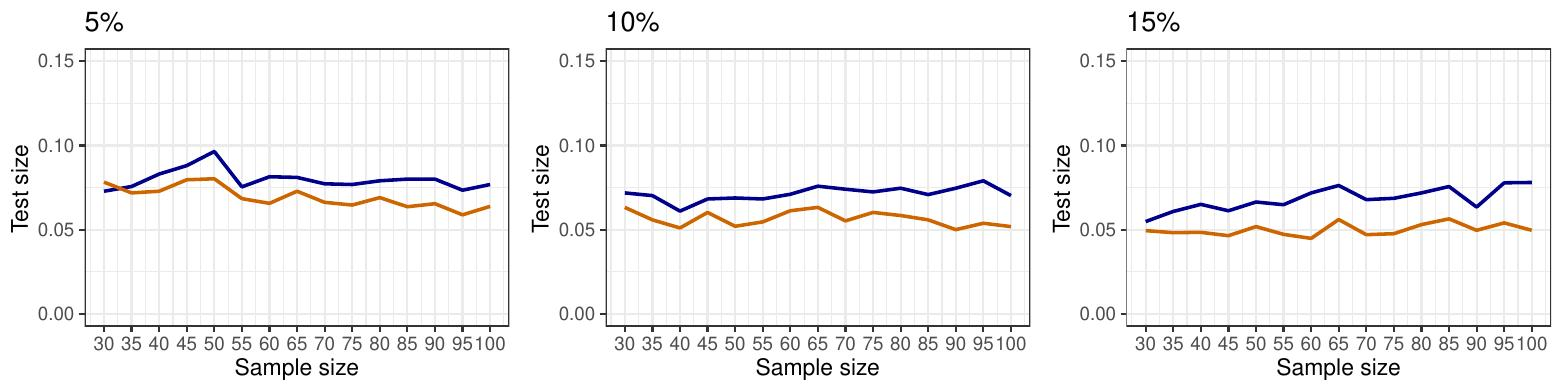}
    \caption{Empirical test \textbf{sizes for 2X3Y case} ($A_n$ in orange, $d^2$ in blue) for \textbf{Clayton copula with parameter 1} and $\mathbf{\chi^2_4}$ \textbf{margins}, \textbf{sample size ranging from 30 to 100}, $N = 5000$ simulations }
    \label{fig:sizes_2X3Y_clayton1_chisq4_n30to100}
\end{figure}

When the underlying distribution is not normal both tests experience slower convergence to the asymptotic distribution, especially in the 2X3Y case. For Clayton copula with parameter 1 and $\chi^2_4$ margins, these observations can be made from figures \ref{fig:sizes_1X2Y_clayton1_chisq4_n30to100}, \ref{fig:sizes_3X2Y_clayton1_chisq4_n30to100} and \ref{fig:sizes_2X3Y_clayton1_chisq4_n30to100}. We see that $A_n$ performs better than Little's $d^2$ in this case as well.

\begin{figure}[ht]
    \centering
    \includegraphics[width=\textwidth]{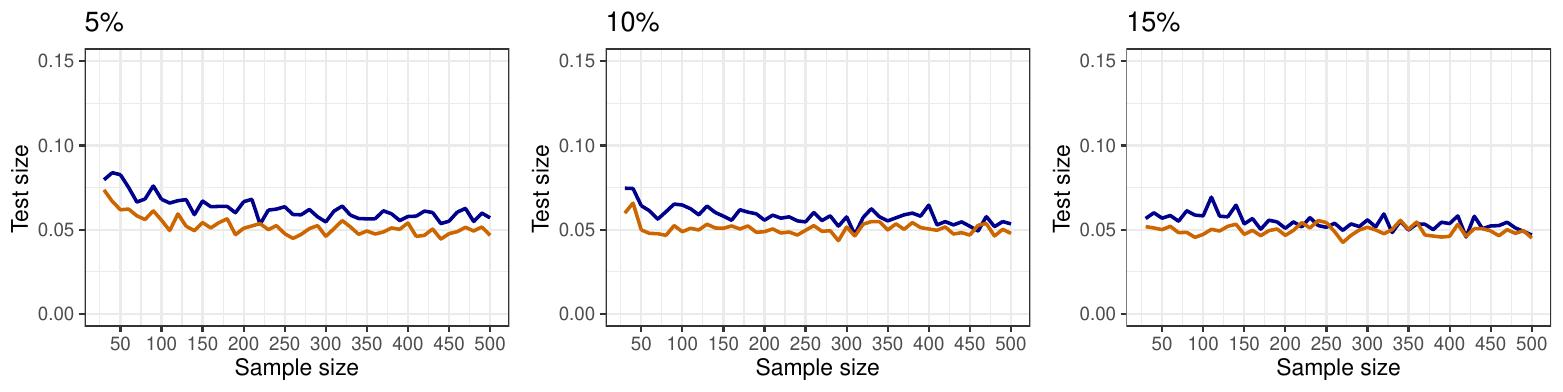}
    \caption{Empirical test \textbf{sizes for 1X2Y case} ($A_n$ in orange, $d^2$ in blue) for \textbf{Clayton copula with parameter 1} and $\mathbf{\mathcal{E}(1)}$ \textbf{margins}, \textbf{sample size ranging from 30 to 500}, $N = 5000$ simulations }
    \label{fig:sizes_1X2Y_clayton1_exp1_n30to500}
\end{figure}

\begin{figure}[ht]
    \centering
    \includegraphics[width=\textwidth]{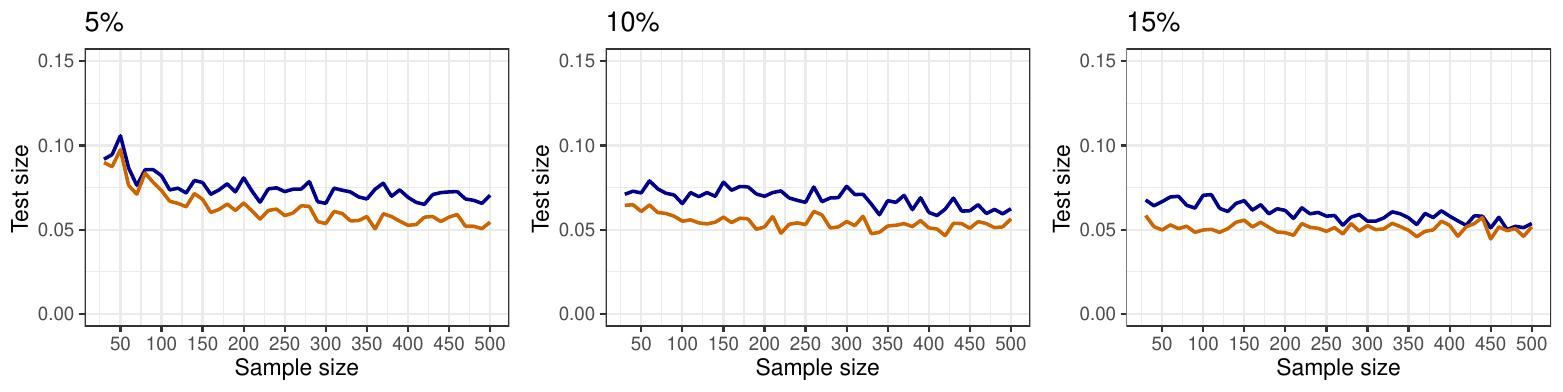}
    \caption{Empirical test \textbf{sizes for 3X2Y case} ($A_n$ in orange, $d^2$ in blue) for \textbf{Clayton copula with parameter 1} and $\mathbf{\mathcal{E}(1)}$ \textbf{margins}, \textbf{sample size ranging from 30 to 500}, $N = 5000$ simulations }
    \label{fig:sizes_3X2Y_clayton1_exp1_n30to500}
\end{figure}

\begin{figure}[ht]
    \centering
    \includegraphics[width=\textwidth]{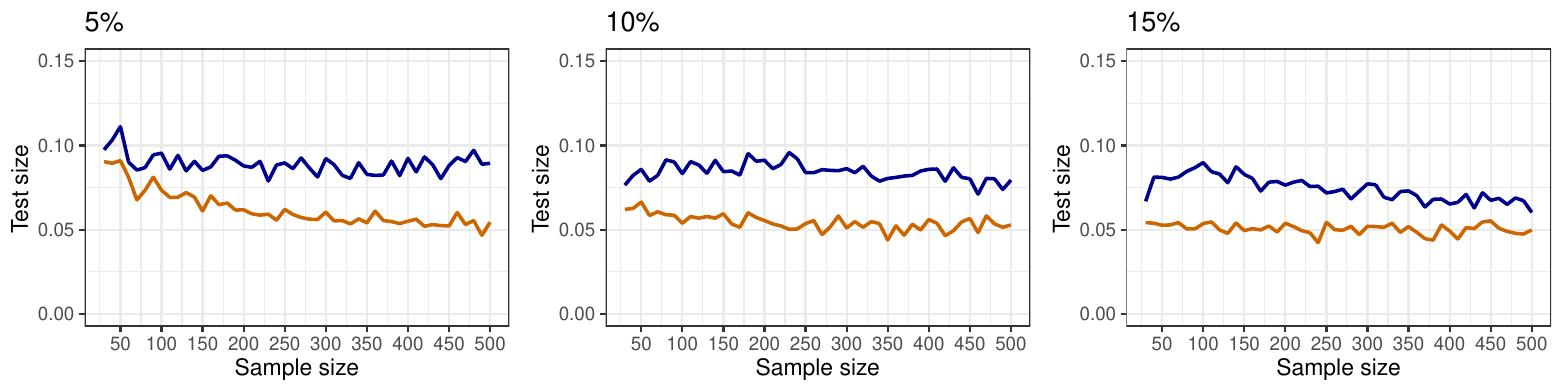}
    \caption{Empirical test \textbf{sizes for 2X3Y case} ($A_n$ in orange, $d^2$ in blue) for \textbf{Clayton copula with parameter 1} and $\mathbf{\mathcal{E}(1)}$ \textbf{margins}, \textbf{sample size ranging from 30 to 500}, $N = 5000$ simulations }
    \label{fig:sizes_2X3Y_clayton1_exp1_n30to500}
\end{figure}

As we can see from figures \ref{fig:sizes_1X2Y_clayton1_exp1_n30to500}, \ref{fig:sizes_3X2Y_clayton1_exp1_n30to500} and \ref{fig:sizes_2X3Y_clayton1_exp1_n30to500}, when the underlying distribution is Clayton copula with parameter 1, but with $\mathcal{E}(1)$ margins, this difference is even more significant. Test based on $d^2$ experiences little to no improvement as sample size increases, even when sample sizes up to 500 were used, unlike with the previous distributions, where sample of 100 units yielded acceptable results. This is especially true for the 2X3Y case (Figure \ref{fig:sizes_2X3Y_clayton1_exp1_n30to500}).

Having these results, we proceed by computing tests' powers using some feasible sample sizes for different distributions ($100$ and $300$ for standard normal distribution, and $300$ for Clayton copula with parameter 1 and $\chi^2_4$ margins).

\subsection{Power comparison}

Under the feasible circumstances, i.e. underlying distribution and sample size that preserve the type I error, we proceed by comparing empirical test powers for statistics $A_n$ and $d^2$.

\begin{figure}[ht]
    \centering
    \includegraphics[width=\textwidth]{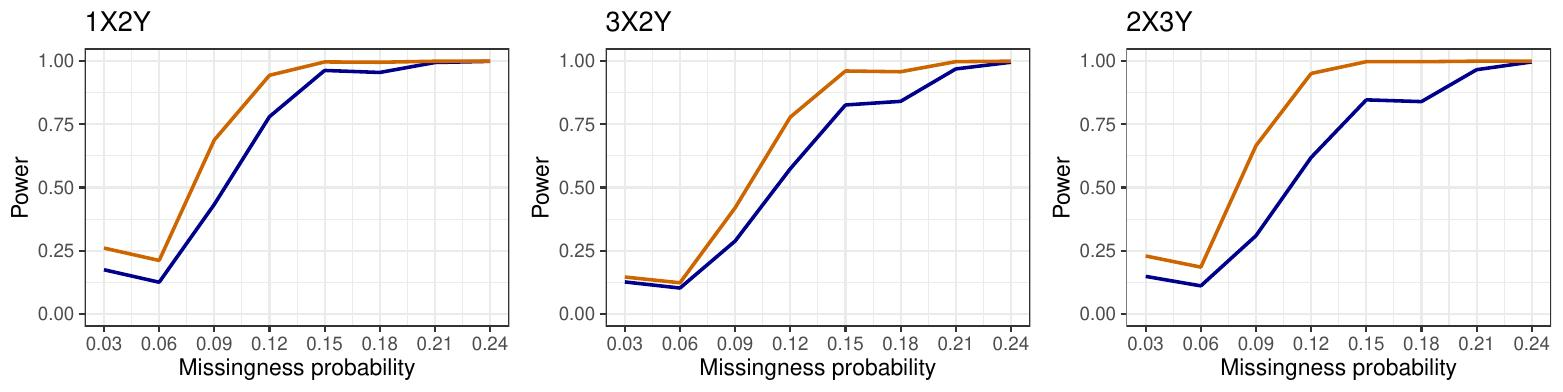}
    \caption{Empirical test \textbf{powers} ($A_n$ in orange, $d^2$ in blue) for \textbf{standard normal} distribution, sample size $\mathbf{n = 100}$, and $N = 5000$ simulations, and \emph{\textbf{MAR 1 to 9}} alternative}
    \label{fig:powers_stdnorm_n100_MAR1to9}
\end{figure}

\begin{figure}[ht]
    \centering
    \includegraphics[width=\textwidth]{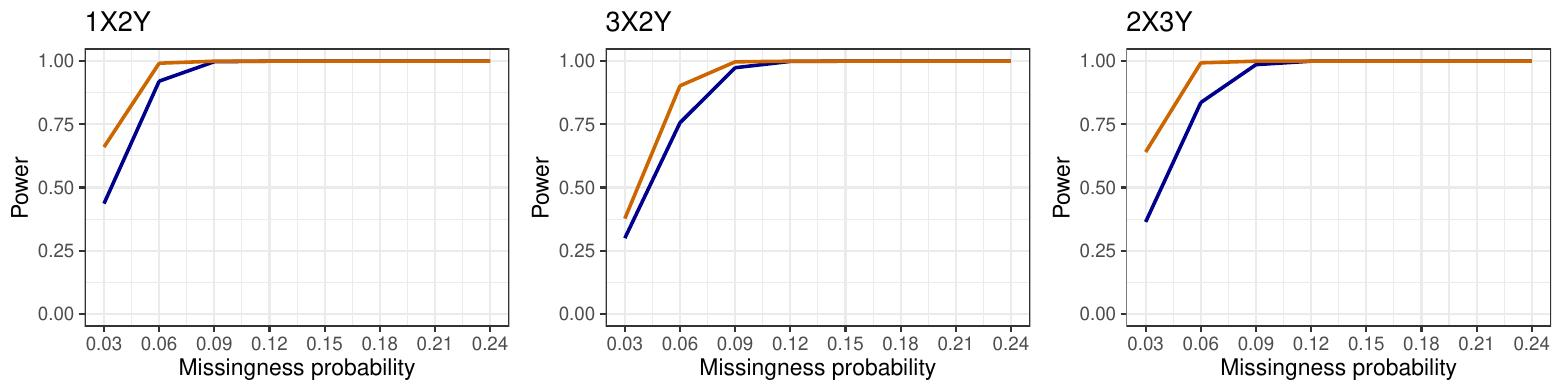}
    \caption{Empirical test \textbf{powers} ($A_n$ in orange, $d^2$ in blue) for \textbf{standard normal} distribution, sample size $\mathbf{n = 300}$, and $N = 5000$ simulations, and \emph{\textbf{MAR 1 to 9}} alternative}
    \label{fig:powers_stdnorm_n300_MAR1to9}
\end{figure}

Since both tests preserve the size for underlying standard normal distribution and sample size of $n = 100$, it is acceptable to compare powers in that case. Figure \ref{fig:powers_stdnorm_n100_MAR1to9} shows that test constructed around $A_n$ is more powerful for all dimensions and patterns under the \emph{MAR 1 to 9} alternative. As expected, for larger sample size all powers increase, but the relationship remains unchanged. This can bee seen in Figure \ref{fig:powers_stdnorm_n300_MAR1to9}.

\begin{figure}[ht]
    \centering
    \includegraphics[width=\textwidth]{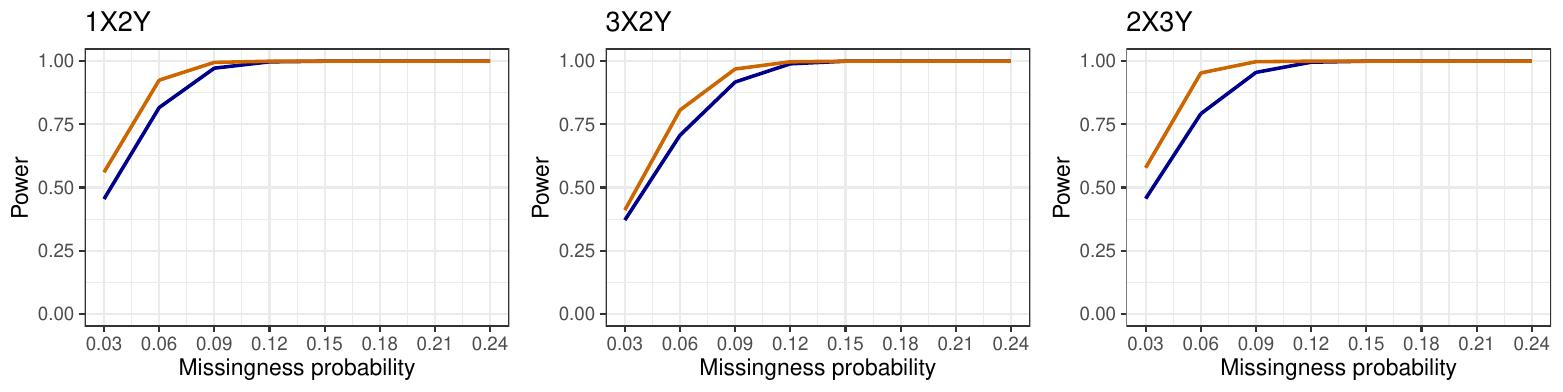}
    \caption{Empirical test \textbf{powers} ($A_n$ in orange, $d^2$ in blue) for \textbf{Clayton copula with parameter 1} and $\mathbf{\chi^2_4}$ \textbf{margins}, sample size $\mathbf{n = 300}$, and $N = 5000$ simulations, and \emph{\textbf{MAR 1 to 9}} alternative}
    \label{fig:powers_clayton1_chisq4_n300_MAR1to9}
\end{figure}

In Figure \ref{fig:powers_clayton1_chisq4_n300_MAR1to9}, we see that for the underlying distribution that is Clayton copula with parameter 1 and $\chi^2_4$ margins, we again have that $A_n$ performs more powerful than $d^2$ under the \emph{MAR 1 to 9}, but the difference seems to be less notable than for the standard normal case. However, one needs to take these powers with caution: Little's test shows tendency to reject null hypothesis, so the powers that we will observe can be considered slightly lower than those obtained, especially for the $2X3Y$ case. That leads us to the same conclusion as for the underlying standard normal distribution.

\begin{figure}[ht]
    \centering
    \includegraphics[width=\textwidth]{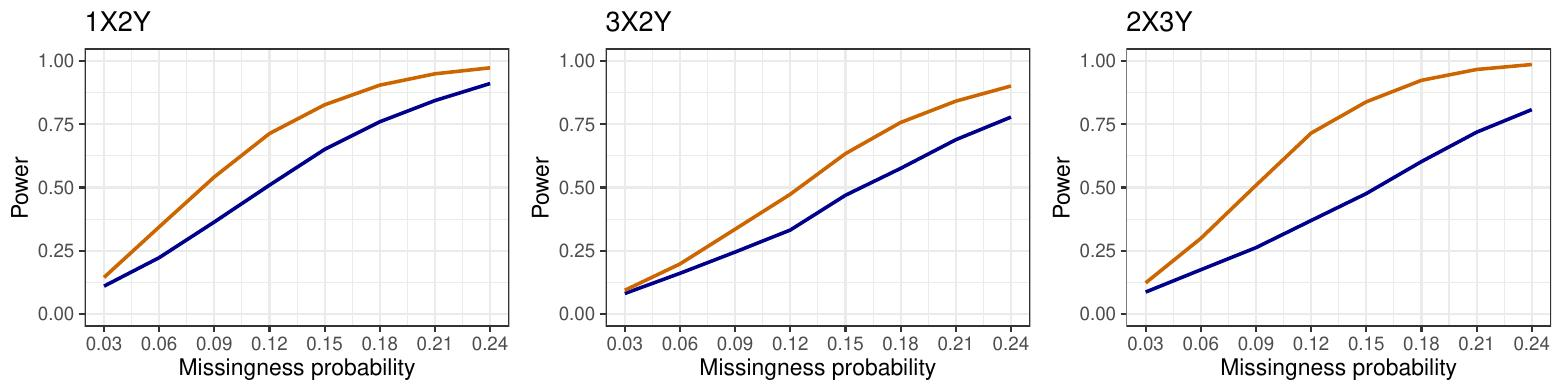}
    \caption{Empirical test \textbf{powers} ($A_n$ in orange, $d^2$ in blue) for \textbf{standard normal} distribution, sample size $\mathbf{n = 100}$, and $N = 5000$ simulations, and \emph{\textbf{MAR rank}} alternative}
    \label{fig:powers_stdnorm_n100_MARrank}
\end{figure}

\begin{figure}[ht]
    \centering
    \includegraphics[width=\textwidth]{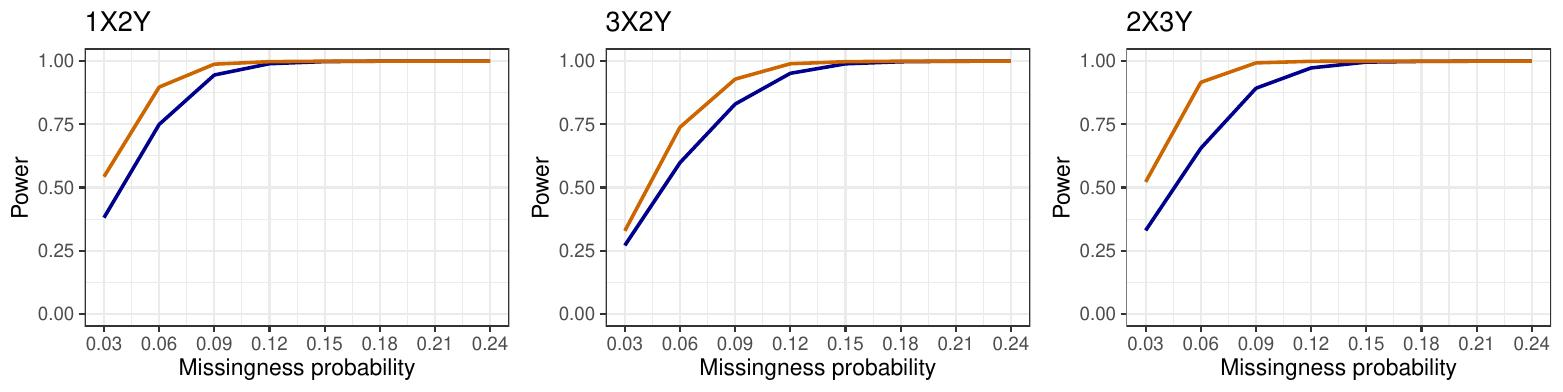}
    \caption{Empirical test \textbf{powers} ($A_n$ in orange, $d^2$ in blue) for \textbf{standard normal} distribution, sample size $\mathbf{n = 300}$, and $N = 5000$ simulations, and \emph{\textbf{MAR rank}} alternative}
    \label{fig:powers_stdnorm_n300_MARrank}
\end{figure}

\begin{figure}[ht]
    \centering
    \includegraphics[width=\textwidth]{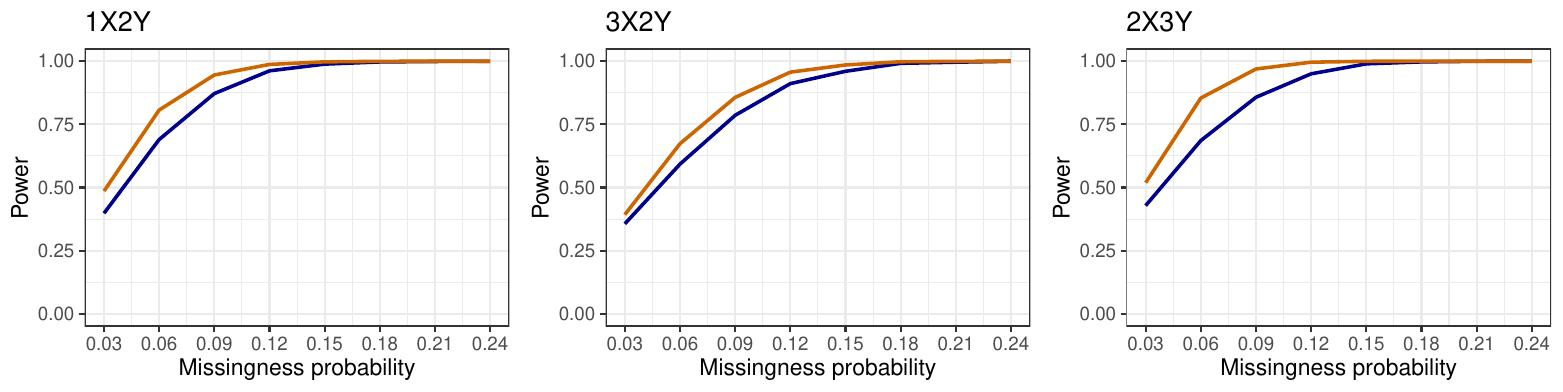}
    \caption{Empirical test \textbf{powers} ($A_n$ in orange, $d^2$ in blue) for \textbf{Clayton copula with parameter 1} and $\mathbf{\chi^2_4}$ \textbf{margins}, sample size $\mathbf{n = 300}$, and $N = 5000$ simulations, and \emph{\textbf{MAR rank}} alternative}
    \label{fig:powers_clayton1_chisq4_n300_MARrank}
\end{figure}

For the \emph{MAR rank} mechanism, we experience that, for the underlying standard normal distribution, empirical powers of both tests increase slower than for the \emph{MAR 1 to 9}, but the relation between tests remains unchanged: $A_n$ is more powerful, which can be seen in Figure \ref{fig:powers_stdnorm_n100_MARrank}. Figure \ref{fig:powers_stdnorm_n300_MARrank} shows expected behavior for increased sample size of $n=300$. Results for the Clayton copula with parameter 1 and $\chi^2_4$ margins can be seen in the Figure \ref{fig:powers_clayton1_chisq4_n300_MARrank}, and they behave as the previous. We note that several other dimensions and underlying distributions were tried. Similarly to the discussion about the sizes, since the observed behavior was identical to this presented here, it is redundant to include them.

\begin{figure}[ht]
    \centering
    \includegraphics[width=\textwidth]{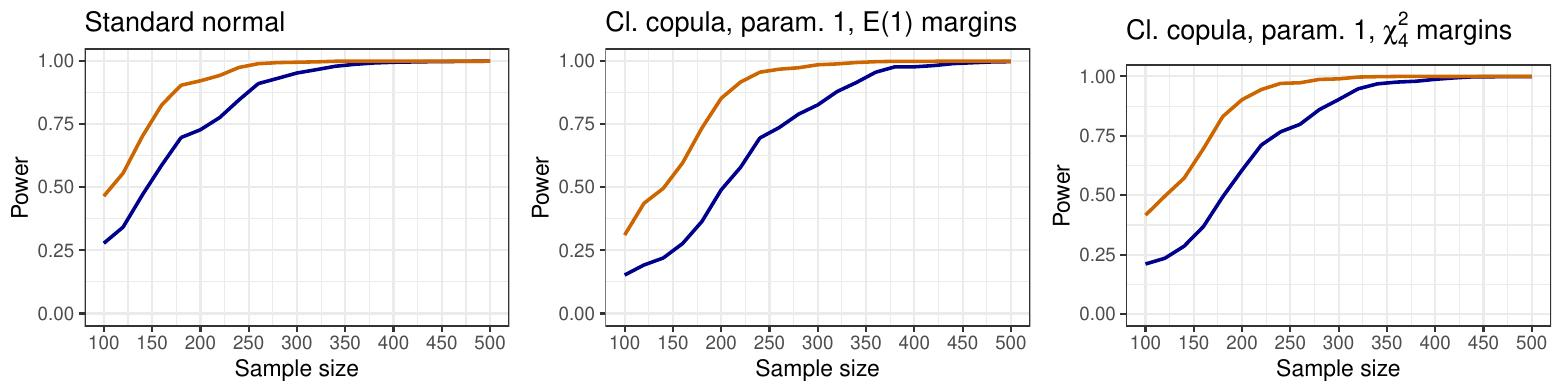}
    \caption{Empirical test \textbf{powers} ($A_n$ in orange, $d^2$ in blue) for \emph{\textbf{MAR mean}} alternative, sample size ranging from 100 to 500, $N = 5000$ simulations}
    \label{fig:powers_1X2Y_n100to500_depth}
\end{figure}

For the \emph{MAR mean} mechanism from the \cite{mozharovskyi2020nonparametric}, the original missingness probabilities needed to be changed. Originally, second column had missingness probability of $0.48$ for those cases where $X$ is larger than its mean value, and $0.24$ for those that have it smaller.  The probabilities for the third column were set to be $0.08$ and $0.7$, respectively. This gives a large number of missing values, that allows both tests to have powers almost equal to 1. To avoid that, we changed the above probabilities to $0.12$, $0.06$, $0.02$, and $0.175$ respectively. We observe the behavior for the standard normal distribution, Clayton copula with parameter 1 $\mathcal{E} (1)$ margins, and Clayton copula with the same parameter and $\chi^2_4$ margins. The sample size is ranged from 100 to 500. One should not that note that, from the previous results, we know that test based on $d^2$ has a tendency to reject the null hypothesis (since sizes are slightly larger than $0.05$). That means that true values of powers for the Little's MCAR test can be considered slightly lower than those observed. Results are presented in Figure \ref{fig:powers_1X2Y_n100to500_depth}. As one can see, results follow the same pattern seen before: our novel test experiences higher powers than the Little's test.

\section{Concluding remarks}

We introduced a novel test for testing the MCAR assumption that is based on the test statistic $A_n$ and we utilized existing theory about $U$-statistics to derive its asymptotic properties. A extensive simulation study was conducted to examine how novel test behaves in the sense of size preservation, speed of convergence to the derived asymptotic $\chi^2$ distribution, and in terms of power. The test has shown better properties both in the preservation of the type I error and in terms of the empirical power of the test, when compared to the well-known Little's MCAR test. The main conclusion is that Little's test is less robust in the sense of underlying distribution, and experiences significant loss in the size preservation when it is not normal. Our novel test also experiences damage when normality assumption is not met, but that damage is much less significant when compared to the other one.

The only visible drawback of the novel test is that it needs at least one complete column to perform. However, situations where all of the columns have missing data are not that common in practice, especially for the big data cases with large number of variables.

\section*{Acknowledgments}

The author would like to thank professor Bojana Milo\v sevi\' c from the University of Belgrade, Faculty of Mathematics, for fruitful discussions throughout creation of this paper.

\subsection*{Financial disclosure}

The work of D. Aleksi\'c is supported by the Ministry of Science, Technological Development and Innovations of the Republic of Serbia (the contract 451-03-47/2023-01/ 200151). 

\subsection*{Conflict of interest}

The author declares no potential conflict of interests.

\newpage
\section*{Figures}

\end{document}